\documentclass[a4paper,11pt]{article}

\usepackage{amsmath}
\usepackage{amssymb}
\usepackage{amsthm}
\usepackage{graphicx}
\usepackage[usenames]{xcolor}
\usepackage[noadjust]{cite}
\usepackage[T1]{fontenc}
\usepackage{todonotes}
\usepackage[utf8]{inputenc}
\usepackage{authblk}
\usepackage[hyphens]{url}
\usepackage[final]{changes}

\usepackage{mathtools}
\usepackage{enumitem}
\usepackage{mathscinet}
\usepackage{bm}

\usepackage{subcaption}
\usepackage[width=136mm]{caption}	

\usepackage[colorlinks=true]{hyperref}
\definecolor{ForestGreen}{rgb}{0.15,0.416,0.18}
\definecolor{EgyptBlue}{rgb}{0.063,0.2,0.65}
\hypersetup{
	colorlinks=true,
	linkcolor=EgyptBlue,         
	citecolor=EgyptBlue,          
	urlcolor=ForestGreen          
}

\usepackage[capitalize,noabbrev]{cleveref}

\newtheorem{theorem}{Theorem}[section]

\newtheorem{lemma}[theorem]{Lemma}
\newtheorem{proposition}[theorem]{Proposition}

\theoremstyle{definition}
\newtheorem{definition}[theorem]{Definition}
\theoremstyle{definition}
\newtheorem{remark}[theorem]{Remark}
\theoremstyle{definition}

\numberwithin{equation}{section}
\numberwithin{table}{section}
\numberwithin{figure}{section}

\def\U{\mathcal{U}}
\def\S{\mathcal{S}}
\def\M{\mathcal{M}}
\def\C{\mathcal{C}}
\def\A{\mathcal{A}}
\def\V{\mathcal{V}}

\def\N{\mathbb{N}}
\def\R{\mathbb{R}}

\def\Anull{\A \setminus \{ \0 \}}
\def\sgn{\operatorname{sgn}}

\newcommand{\eps}{\ensuremath{\varepsilon}}
\newcommand{\scc}{\operatorname{sc}}
\newcommand{\Vmax}{\ensuremath{V_{\mathrm{max}}}}

\newcommand{\Lip}[1]{\ensuremath{\operatorname{Lip}({#1})}}

\newcommand{\naturals}{\ensuremath{\mathbb{N}}}

\newcommand{\complex}{\ensuremath{\mathbb{C}}}

\newcommand{\norm}[1]{\ensuremath{\Vert#1\Vert}}
\newcommand{\setof}[1]{\ensuremath{\left\{#1\right\}}}

\newcommand{\intd}{\ensuremath{\mathrm{d}}}

\newcommand{\phasespace}{\ensuremath{\mathcal{X}}}
\newcommand{\0}{\ensuremath{\mathbf{0}}}
\newcommand{\semiflow}{\ensuremath{\mathcal{F}}}

\makeatletter
\newcommand{\leqnomode}{\tagsleft@true\let\veqno\@@leqno}
\newcommand{\reqnomode}{\tagsleft@false\let\veqno\@@eqno}
\makeatother

\makeatletter
\newcommand{\subjclass}[2][2020]{%
  \let\@oldtitle\@title%
  \gdef\@title{\@oldtitle\footnotetext{#1 \emph{MSC Classification:} #2}}%
}
\newcommand{\keywords}[1]{%
  \let\@@oldtitle\@title%
  \gdef\@title{\@@oldtitle\footnotetext{\emph{Keywords:} #1}}%
}
\makeatother

\title{Morse decomposition of scalar differential equations with state-dependent delay}
\author[1,2,3]{Ferenc A. Bartha\footnote{Corresponding author: barfer@math.u-szeged.hu.}}
\author[1,2,3]{Ábel Garab}
\author[1,2,3]{Tibor Krisztin}

\affil[1]{Bolyai Institute, University of Szeged,\par Aradi vértanúk tere 1, Szeged, H–6720, Hungary}

\affil[2]{HUN-REN--SZTE Analysis and Applications Research Group, \par Bolyai Institute,
University of Szeged, Szeged, Hungary}

\affil[3]{National Laboratory for Health Security, University of Szeged, Szeged, Hungary}



\date{}

\begin{document}

%%%%%%%%%%%%%%%%%%%%%%%%%%%%%%%%%%%%%%%%%%%%%%%%%%%%%%%%%%%%%%%%%%%%%%%%
%%%%%%%%%%%%%%%%%%%%%%%%%%%%%%%%%%%%%%%%%%%%%%%%%%%%%%%%%%%%%%%%%%%%%%%%
\maketitle
\begin{abstract}
\noindent We consider state-dependent delay differential equations of the form 
	\begin{equation*}
	    \dot{x}(t) = f(x(t), x(t - r(x_t))),
	\end{equation*}
	where $f$ is continuously differentiable and fulfills a negative feedback condition in the delayed term. Under suitable conditions on $r$ and $f$, we construct a Morse decomposition of the global attractor, giving some insight into the global dynamics. The Morse sets in the decomposition are closely related to the level sets of an integer valued Lyapunov function that counts the number of sign changes along solutions on intervals of length of the delay. This generalizes former results for constant delay. %by Mallet-Paret [\emph{J. Differential Equations}, \textbf{72} (1988), 270--315] and Polner [\emph{Nonlinear Analysis}, (2002) 377--397]. 
    We also give two major types of state-dependent delays for which our results apply.

\end{abstract}

\section{Introduction}
\label{sec:introduction}

Let us consider delay differential equations of the form 
\begin{equation}\label{eq:intro}
    \dot{x}(t) = f(x(t), x(t - r)),
\end{equation}
where $f$ is continuously differentiable and fulfills  a negative feedback condition in the delayed term.

Such equations have been studied in great detail when the delay $r$ is constant 
\cite{arino, fiedler, krisztin-walther, krisztin-walther-wu, mallet-paret, mallet-sell, mallet-sell2, mischaikow}. 
On the other hand, if $r$ depends on the solution $x(t)$, we speak of a state-dependent delay. 
Such equations arise naturally in e.g., models of automatic position control, population dynamics, 
mill-control, white blood cells, see \cite{hartung-krisztin-walther-wu} and the references therein. 
These equations do not fit in the classical theory of functional-differential equations \cite{hale-lunel}. 
The usual smoothness properties are violated due to the state-dependent delay. 
%\added{Even though several important results have been achieved, e.g.\ generic structure of the global attractor \cite{krisztin-arino}, existence of periodic solutions, phase--plane analysis \cite{mallet-nussbaum1}, geometry of the solution manifold \cite{krisztin-walther-almost-graph,walther};}  
%numerous questions remain open for these systems \cite{krisztin-walther-smoothness}.

\added{
A basic theory has been developed by Walther for well-posedness of initial value problems within the dynamical systems framework on solution manifolds in spaces of differentiable functions \cite{waltherC1,hartung-krisztin-walther-wu}. Although several problems could be solved within this framework, it has limitations. As the solution operators are in general only $C^1$-smooth, higher smoothness of local invariant manifolds is not available \cite{krisztinunstable,krisztincenter,quesmi-walther,stumpfcenter,krisztin-walther-smoothness}. There is some information on the geometry of solution manifolds \cite{krisztin-walther-almost-graph,walther}. 
Considerable progress has been made in understanding some prototype examples. We mention only a few here. Mallet-Paret and Nussbaum obtained boundary layer phenomena in singularly perturbed systems \cite{mallet-nussbaum1,mallet-nussbaum2,mallet-nussbaum3,mallet-nussbaum4}. 
Further recent work extended global methods to find periodic orbits (including hyperbolicity) to equations with state-dependent delay \cite{magal-arino,waltherPositionControl,walther-periodic}, identified (parts of) global attractors \cite{krisztin-arino,stumpfcenterunst}, established multiple periodic orbits \cite{kennedy-multiple}, obtained local and global Hopf bifurcation \cite{eichmann,hu-wu}, and proved complicated behavior of solutions \cite{kennedy-chaos,laniwayda-walther-shilnikov}. 
In spite of the above mentioned works and several other important results, numerous fundamental questions remain open.
}

In particular, Kennedy recently proved a Poincar\'e--Bendixson-type result for the case of monotone negative feedback \cite{kennedy}. This is a scalar, state-dependent delay analog of the result obtained by Mallet-Paret and Sell in their pioneering work \cite{mallet-sell2} for cyclic systems with constant delay. In the general case, i.e.~when monotonicity in the second argument of $f$ is not assumed, one does not expect such a strong result to hold. However, in case of constant delay and under some additional assumptions on the right-hand side, Mallet-Paret \cite{mallet-paret} and later Polner \cite{polner} established the existence of a so-called Morse decomposition of the global attractor -- based on an integer-valued Lyapunov function --  for negative and positive feedback, respectively. A Morse decomposition is a finite
collection of pairwise disjoint, compact invariant subsets of the global attractor,  which are ordered in the sense that -- roughly speaking -- the dynamics on the
attractor and outside these sets is gradient-like \cite{conley}. Note that this does not necessarily imply simplicity: inside a Morse component the system may exhibit complicated behavior, even chaos -- see the papers \cite{laniwayda-walther, laniwayda-walther_errata,peters,siegberg} in this direction.

The aim of this paper is to give sufficient conditions in terms of $f$ and $r$ for the existence of a Morse decomposition analogous to those in \cite{mallet-paret, polner, garab-poetzsche, garab}. 

Throughout this paper, we use the following notation. 
Let $K > 0$ be fixed. If $t_1 < t_2$ and $u \colon [t_1 - K, t_2) \to \R $ is continuous, 
then, for $t \in [t_1, t_2)$, $u_t \in C([-K, 0], \R )$ is defined by $u_t(s) = u(t + s)$, 
$-K \leq s \leq 0$. 

In our setting, the delay $r$ is defined by the $K$-long segment $x_t$ of the solution $x$, that is,  
\begin{equation}
\label{eq:dde}
    \dot{x}(t) = f(x(t), x(t - r(x_t)). 
\end{equation}
%In the rest of the paper we use the notation $r = r(x, t)$ interchangeably with $r = r(x_t)$ if it makes the reasoning simpler to follow. 

Some relevant ideas for \eqref{eq:dde} have already been introduced in 
\cite{krisztin-arino, mallet-nussbaum1, mallet-nussbaum-para}, however, 
the particular form $r = r(x(t))$ plays an important role in them.
Our result is new for the case $r = r(x(t))$ as well. 

Let $M > 0$ be fixed. 
We make the following assumptions on $f$:  
\begin{enumerate}[label=$(\mathrm{H}_\arabic*)$,  leftmargin=*]
\item   $\displaystyle f\in C^1(\R ^2, \R )$ and $\displaystyle L_0\coloneqq \sup_{x, y \in [-M, M]} |f(x,y)| ,$     \label{eq:f:bounded} \\ 
\item	%$\displaystyle \frac{\partial f(x,y)}{\partial y}(0,0)<0 \quad \text{and} \quad y f(0, y) < 0 \quad \mbox{for all} ~ y \neq 0, $
$\displaystyle
y f(0, y) < 0$ for all $y \neq 0, \quad \text{and}
\quad \frac{\partial f(x,y)}{\partial y}(0,0)<0  $,
\label{eq:f:negative} 
\item $x f(x, y) < 0$ for all $(x, y) \in \setof{(u, v) \in \R ^2 \colon |u| \geq \max \! \setof{M, |v|}}$. 
     \label{eq:f:dissipative}
\end{enumerate}
Hypothesis \ref{eq:f:negative} is the negative feedback condition in the second variable. 
That, together with the continuity of $f$ implies
\begin{equation}
\label{eq:f:zero}
	f(0, 0) = 0. 
\end{equation}
\added{Note that the second part of \ref{eq:f:negative} assumes solely local monotonicity at the origin in the second variable.} 
Property \ref{eq:f:dissipative} is referred to as the dissipativity condition. 
\added{As an example, $f(x,y)=-ax-g(y)$ with $a>0$ and a bounded continuously differentiable real function $g$ with $g'(0)>0$, and  $yg(y)>0$ for all $y\ne 0$ satisfies the conditions \ref{eq:f:bounded}, \ref{eq:f:negative}, \ref{eq:f:dissipative}.}

Solutions will be defined in Section \ref{sec:solutions-delays}. In the sequel, we shall only consider solutions of \eqref{eq:dde} such that $x(t) \in (-M, M)$ 
for all $t$ in the domain of $x$. 
We adopt the way of defining the phase space from \cite{krisztin-arino}.
First, consider the Banach space $\C = C([-K, 0], \R )$ 
with the norm $\norm{\varphi} = \max_{s \in [-K, 0]} |\varphi(s)|$ 
for $\varphi \in \C$. We will also use the $C^1$-norm: for $\psi\in C^1([a,b],\R )$, let $\|\psi\|_{C^1([a,b],\R )} \coloneqq \sup_{t\in [a,b]} |\psi(t)| + \sup_{t\in [a,b]} |\dot{\psi}(t)|$. For brevity, we introduce the notation $\|\cdot\|_1$ for $\|\cdot\|_{C^1([-K,0],\R )}$. For a map $h\colon D\to \R$, where $D$ is a subset of a Banach space $E$ with norm $|\cdot|_E$, let us define
\[\Lip{h}\coloneqq \sup_{x,y\in D, ~ x \neq y}\frac{|h(x)-h(y)|}{|x-y|_E}\leq \infty.\] 

We choose 
\begin{equation*}
	\phasespace = \setof{\varphi \in \C \colon \norm{\varphi} < M \text{ and } \Lip{\varphi}\leq L_0},
\end{equation*}
as the \textit{phase space} for \eqref{eq:dde}.

Note that by the Arzel\`a--Ascoli theorem, the closure $\overline{\phasespace}$ is a compact subset of $\C$. 

Furthermore, we assume throughout the paper that the state-dependent delay $r \colon \overline{\phasespace} \to \R $ satisfies 
\begin{enumerate}[label=$(\mathrm{H}_\arabic*)$, start=4, leftmargin=*]
\item	$r(\varphi) \in (0, K] \quad \mbox{for all} ~ \varphi \in \phasespace, $
     \label{eq:delay:bounded} 
\item $r$ is globally Lipschitz continuous with $\Lip{r}$, 
     \label{eq:delay:lip}
\item
     $t\mapsto \eta_x(t)\coloneqq t-r(x_t)$  is strictly increasing for any solution $x$.
     \label{eta:increasing}
\end{enumerate}

\added{We will frequently use the notation $\eta^n_x(t)$ ($n\in \N_0)$ for the $n$-th iterate of the function $\eta_x$, that is, $\eta^0_x(t)\coloneqq t$ and $\eta^{m+1}_x(t)=\eta_x(\eta^m_x(t))$, $m\geq 0$. Note that if we had constant delay $r(x_t) \equiv \tau$, then $\eta^n_x(t)$ would simplify to $t-n\tau$.} 

Hypotheses \ref{eq:f:bounded}--\ref{eta:increasing} can be verified for a wide class of equations -- see Section~\ref{sec:examples} for two very common types of state-dependent delays and references \cite{hu-wu,magal-arino,arino-hadeler-hbid} for further examples.

Our main result is that a Morse decomposition of the global attractor exists, provided hypotheses \ref{eq:f:bounded}--\ref{eta:increasing} hold.
\medskip

The paper is organized as follows. 
In \cref{sec:solutions-delays}, we establish the existence of the global attractor and introduce a related linear, nonautonomous equation with variable delay that is closely related to \eqref{eq:dde}.

As already indicated, our result relies on an integer-valued Lyapunov function that essentially counts the number of sign changes of solution $x$ on intervals $[t-r(x_t),t]$. The idea of this concept originates from My\v{s}kis \cite{myshkis} and was refined by Mallet-Paret and Sell \cite{mallet-sell} for constant delay. We shall use a state-dependent analogue of this Lyapunov function, and extensively use its properties. Many of these were already established by Krisztin and Arino \cite{krisztin-arino}, however, some of their statements had to be extended to fit our purposes. This preparatory work is included in \cref{sec:signchange}. 

\Cref{sec:morse} is dedicated to the proof of our main result, \cref{thm:morse}, in which a Morse decomposition of the global attractor is given. The results are analogous to those with constant delay \cite{polner, mallet-paret,garab}. The technical difficulties arising due to the state-dependency of the delay are nontrivial. The main difference in the result compared to the constant delay case is that at this level of generality of the delay we do not expect the discrete Lyapunov function to be bounded on the global attractor. Yet, even in this case, there exists a Morse decomposition with finitely many Morse sets, one of which contains the segments of
all entire solutions with oscillation frequencies larger than a prescribed frequency. This set may contain all possible superhigh-frequency oscillations, i.e. the ones that change sign infinitely often on certain bounded intervals.

In \cref{sec:sharper:results} we give condition \ref{iterated_zeros_imply_fullzero} guaranteeing boundedness of the discrete Lyapunov function on the global attractor.

Finally, in \cref{sec:examples}  we restrict our attention to two classes of state-dependent delays. We show that these fulfill hypotheses \ref{eq:delay:bounded}--\ref{eta:increasing}, so our results apply. Moreover, we demonstrate that for these two classes of delays, our integer valued Lyapunov function is bounded on the global attractor. This allows for a sharper result (\cref{thm:morse:spec:delays}), completely analogous to the constant delay case.

%%%%%%%%%%%%%%%%%%%%%%%%%%%%%%%%%%%%%%%%%%%%%%%%%%%%%%%%%%%%%%%%%%%%%%%%
%%%%%%%%%%%%%%%%%%%%%%%%%%%%%%%%%%%%%%%%%%%%%%%%%%%%%%%%%%%%%%%%%%%%%%%%
\section{Preliminaries}
\label{sec:solutions-delays}

\subsection{Solutions and the global attractor}
We call the continuous function $x \colon [-K, t_\varphi) \to (-M, M)$ a \textit{solution of} \eqref{eq:dde} starting 
from the initial function $\varphi \in \phasespace$, if it is differentiable and satisfies \eqref{eq:dde} 
on $(0, t_\varphi)$ and $x_0 = \varphi$.  The Lipschitz continuity of $r$ and of $\varphi \in \phasespace$ imply that if a solution $x$ starting 
from a $\varphi \in \phasespace$ exists, then it is unique, hence, we denote it by $x^{\varphi}$. %ide jo az uj paragraph hogy x^fi kitunjon a szovegbol

For a solution $x\colon I\to \R$ of \eqref{eq:dde}, where $I$ is an interval, we will often use the notation $\eta_x(t)\coloneqq t-r(x_t)$ to denote the delayed argument.

We have continuous dependency on initial conditions: for $\varphi \in \phasespace$, $t \in [0,t_\varphi]$, and $\eps > 0$, 
there exists $\delta > 0$, such that for all $\psi\in \phasespace$ with $\|\varphi - \psi\|<\delta$ the solution $x^\psi$ exists on $[-K,t_\psi)$ with $t_\psi>t$, and $|x^\psi(s) - x^\varphi(s)|<\eps$ for all $s\in [-K,t]$.

Neither uniqueness nor existence is guaranteed backwards. Nevertheless, we will often work with 
solutions of \eqref{eq:dde} defined on the whole real line: a differentiable function $x\colon \R  \to (-M,M)$ satisfying \eqref{eq:dde} 
for all $t\in \R $ is called an \textit{entire solution}. For an entire solution $x \colon \R  \to (-M,M)$ {\em through} $\varphi$, i.e.\ with $x_0 = \varphi$, we also use the notation $x^\varphi$.

The following result ensures the existence of solutions on $[-K, \infty)$ and leads to the existence of a unique connected attractor of the semiflow generated by \eqref{eq:dde}.

\begin{proposition}
\label{prop:semiflow}
	For all $\varphi \in \phasespace$, 
	\eqref{eq:dde} has a unique solution  
	$x^\varphi \colon [-K, \infty) \to \R  $. Moreover,
	\begin{equation*}
		\semiflow \colon [0, \infty) \times \phasespace \ni 
		 (t, \varphi) \mapsto x^\varphi_t \in \phasespace
	\end{equation*}
	is a continuous semiflow.
\end{proposition}

\begin{proof}
    First, for $\varphi\in \phasespace$, let 
    \begin{equation}
    \label{eq:F-rhs}
        F(\varphi) = f(\varphi(0), \varphi(-r(\varphi))).
    \end{equation}
    Then, for $\varphi,\psi \in \phasespace$, we have
	\begin{equation*}
	\begin{split}
		\left| F(\varphi) - F(\psi) \right| &\leq 
		\Lip{f} \left| \varphi(0) - \psi(0) \right| + \Lip{f} \left| \varphi(-r(\varphi)) - \psi(-r(\psi)) \right| \\
		&\leq \Lip{f} \norm{\varphi - \psi} + 
        \Lip{f} \left| \varphi(-r(\varphi)) - \psi(-r(\varphi)) \right| \\
        &\mathop{\hphantom{=}}{}+ \Lip{f} \left| \psi(-r(\varphi)) - \psi(-r(\psi)) \right| \\
		&\leq 2  \Lip{f} \norm{\varphi - \psi} +
		\Lip{f} \Lip{\psi} \Lip{r} \norm{\varphi - \psi} \\
		&\leq \Lip{f} (2 + L_0 \Lip{r})  \norm{\varphi - \psi},
	\end{split}
	\end{equation*}
    where $\Lip{f}$ is taken on $[-M, M]^2$ for the continuously differentiable $f$. 
	The Lipschitz property of $F$ implies that for any $\varphi\in \phasespace$, equation \eqref{eq:dde} has a unique maximal solution $x\colon [-K, \alpha)$, $\alpha\in (0,\infty]$, such that $x_0=\varphi$. For details, see \cite{diekmann,hale-lunel}. 
    
    Furthermore, continuity of $f$ and \ref{eq:f:dissipative} imply that there exists a $\delta\in (0,M-\|\varphi\|)$ such that 
	\begin{equation}\label{delta:dissip}
	xf(x,y)<0 \quad \mbox{holds for all}\quad  (x,y)\in \setof{(u, v) \in \R ^2 \colon |u| \geq \max \! \setof{M-\delta, |v|}}.
	\end{equation}
	We claim that $\|x_t\|\leq M-\delta$ for all $t\in [0,\alpha)$. If this is not the case, then there exists a minimal $t_0\in (0,\alpha)$, such that $|x(t_0)|=M-\delta$. For simplicity, assume that $x(t_0)>0$ (the negative case is similar). Then the minimality of $t_0$ implies that $\dot{x}(t_0)\geq 0$. On the other hand, also by the definition of $t_0$, $|x(t_0-r(x_{t_0}))|<M-\delta=x(t_0)$ holds, which together with \eqref{delta:dissip} yield $\dot{x}(t_0)<0$, a contradiction.
	
	From this one readily obtains by a standard continuation argument that $\alpha=\infty$, which completes the proof.
\end{proof}
Note that, by \eqref{eq:f:zero} and \eqref{eq:F-rhs}, we have that  
$F(\0) = \0$, where $\0 \in \phasespace$ is 
the identically zero function. That is $\0$ is an equilibrium point of $\semiflow$. 

\begin{remark}
If a differentiable function $\chi \colon [-K, \infty) \to \R$ 
satisfies \eqref{eq:dde}, one can show 
-- using the dissipativity condition \ref{eq:f:dissipative} -- 
that eventually it is attracted by $(-M, M)$, 
namely, there is a $t^*$ such that $\chi(t) \in (-M, M)$ for all $t > t^*$.
\end{remark}

Now we turn our attention to the global attractor of $\semiflow$.
\begin{definition}
\label{def:attractor}
	A set $\A  \subset \phasespace$ is called a global attractor if it is compact,  invariant, i.e.\ $\semiflow(t, \A ) = \A$ for all $t \geq 0$, and attracts all subsets of $\phasespace$, that is, for all $B \subset \phasespace$ and for all open $U \supset \A$, there exists $t_0 \geq 0$ such that $\semiflow([t_0, \infty), B) \subset U$.
\end{definition}

The next proposition establishes the existence and uniqueness of the global attractor, and gives its standard characterization.

\begin{proposition}\label{attraktorszakaszok} There exists a unique connected global attractor $\A$ of \eqref{eq:dde} and it allows the following characterization
\[\A = \{ \varphi \in \phasespace : \text{there exists an entire solution }
    x \text{ of  \eqref{eq:dde}, through } \varphi \}.\]
\end{proposition}
\begin{proof}
Let $B\subset \phasespace$ be an arbitrary bounded set. Since $\phasespace$ is precompact, so is $\semiflow(t,B)$ for all $t\geq 0$ and therefore $\semiflow(t,\cdot)$ is completely continuous and also trivially point dissipative, see \cite[p.~16]{hale}.

%there exists $\delta>0$ such that $\|\psi\|< M-\delta$ for all $\phi\in B$. An argument similar to the one in the proof of \Cref{prop:semiflow} shows that $\semiflow(t,B)$ is a subset of the compact set $\{\psi\in \phasespace : \|\psi\|\leq M-\delta\}$. In the proof of \cref{prop:semiflow} we have seen that $\semiflow$ is point-dissipative, thus

Thus the existence of a connected global attractor follows immediately from  \cref{prop:semiflow} and \cite[Theorem 3.4.8]{hale}, whereas \cite[Lemma~2.18\,(d)]{raugel} gives the above characterization of it.
\end{proof}

\added{Observe that the above characterization also implies that $\A \subset C^1([-K,0],\R)$.}
%%%%%%%%%%%%%%%%%%%%%%%%%%%%%%%%%%%%%%%%%%%%%%%%%%%%%%%%%%%%%%%%%%%%%%%%

\subsection{Related nonautonomous linear equations}\label{subsec:trafo}
As already mentioned, our main result, the Morse decomposition of the global attractor, is based on an integer valued Lyapunov function that is closely related to the number of sign changes of a solution $x$ on intervals of the form $[t-r(x_t),t]$. 
 
In this subsection we show that solutions of \eqref{eq:dde} solve also a related nonautonomous linear equation of the form \eqref{felbontas} that can be transformed into the form \eqref{y-felbontas}, where the coefficient function $c$ has constant negative sign. The point of this transformation lays in the fact that the original solution $x$ and its transformation $y$ have the same signs everywhere, and the sign changes of the latter are much easier to keep track of.

\begin{lemma}
\label{felb}
Assume that hypotheses \ref{eq:f:bounded}--\ref{eta:increasing} hold. Then, there exist positive constants $\tau_0,a_0,b_0,b_1,c_0$ and $c_1$, such that the following statements hold for any entire solution $x$ of \eqref{eq:dde}. 
\begin{enumerate}[label=(\roman*)]
\item There exist continuous functions $a,b$ such that $x$ also solves the equation 
\begin{equation}
\label{felbontas}
\dot{x}(t) = a(t)x(t) + b(t)x(\eta_x(t))
\end{equation}
and
\begin{equation}\label{ab-bounds}
  \tau_0\leq r(x_t),\qquad  |a(t)|\leq a_0, \quad \text{and}\quad -b_1 \leq b(t) \leq -b_0<0
\end{equation}
are satisfied for all $t\in \R $.

\item Let $t_0\in \R$ be arbitrary,
\begin{equation}\label{eq:x_to_y_transform}
 y(t) = \exp \left(- \int^t_{t_0} a(s) \,\intd s\right) x(t) \quad \text{and} \quad c(t)=b(t)\exp \left(- \int^t_{\eta_x(t)} a(s) \,\intd s\right).
 \end{equation}
Then, $\sgn y(t)=\sgn x(t)$ holds for all  $t\in \R $, $y$ solves the equation
\begin{equation}
\label{y-felbontas}
\dot{y}(t) = c(t) y( \eta_x (t) ) \qquad \mbox{for} ~ t \in \R , \end{equation}
and 
\begin{equation}\label{c-bounds}
-c_1\leq c(t)\leq -c_0<0
\end{equation}
hold for all $t\in \R$.
\end{enumerate}
\end{lemma}

\begin{proof}
Let $x$ be a fixed entire solution. For all $t\in \R$, let $a (t) = \int^1_0 D_1 f (sx(t),x(\eta(t))) \,\intd s$ and $b (t) = \int^1_0 D_2 f (0,sx(\eta(t))) \,\intd s$, where $D_i$ denotes the differentiation w.r.t.\ the $i$-th argument, and $\eta(t)\coloneqq \eta_x(t)=t-r(x_t)$. Note that 
\[b(t)=\begin{dcases}
    D_2 f(0,0), &\text{if } x(\eta(t))=0,\\
    \frac{f(0,x(\eta(t))}{x(\eta(t))}, &\text{otherwise}.
\end{dcases}\]

\noindent
As $x(t)$ is a solution of (\ref{eq:dde}), we have
\begin{equation*}
\dot{x}(t) = f(x(t),x(\eta(t))) = a(t)x(t) + b(t)x(\eta(t))\quad \text{for all }t\in \R.
\end{equation*}
By continuous differentiability of $f$, the maps $a$ and $b$ are continuous.

Since $\A$ is compact, it follows from the characterization in \cref{attraktorszakaszok} that $$|x(t)|\leq \max_{\varphi \in \A} \|\varphi\|<M$$ holds for all entire solutions $x$ and $t\in \R$. In light of this, the feedback condition \ref{eq:f:negative}, and the boundedness and continuity  of the delay \ref{eq:delay:bounded}--\ref{eq:delay:lip}, we readily obtain the existence of positive constants $\tau_0,a_0,b_0$ and $b_1$ satisfying \eqref{ab-bounds} for any entire solution.

Now, consider $y$ and $c$ as in \eqref{eq:x_to_y_transform}. Then, $\sgn y(t) = \sgn x(t)$ and
\begin{align*}
\dot{y}(t) &= \exp \left(- \int^t_{t_0} a(s) \,\intd s\right) \dot{x}(t) - a(t) \exp \left(- \int^t_{t_0} a(s) \,\intd s\right) x(t)\\
&= \exp \left(- \int^t_{t_0} a(s) \,\intd s\right) \left(a(t)x(t) + b(t)x(\eta(t))\right) - a(t) \exp \left(- \int^t_{t_0} a(s) \,\intd s\right) x(t)\\
&= b(t) \exp \left(- \int^t_{t_0} a(s) \,\intd s\right) x(\eta(t)) = 
b(t) \exp \left(- \int_{\eta(t)}^{t} a(s) \,\intd s - \int^{\eta(t)}_{t_0} a(s) \,\intd s\right)x(\eta(t))\\
&= c(t) y(\eta(t)),\qquad  t\in \R.
\end{align*}
Finally, inequalities \eqref{c-bounds} are satisfied if we set
$c_0= b_0 e^{-Ka_0}$ and $ c_1=b_1 e^{Ka_0}.$
\end{proof}

An important consequence of the above lemma is that, in spite of the lack of backwards-uniqueness, the trivial solution is unique backwards in time. %Another straightforward corollary is obtained for entire solutions.

% \begin{corollary}
% \label{felb}
% If $x \colon \R  \to (-M,M)$ is an entire solution of \eqref{eq:dde}, then there exist continuous functions $c, d, \beta, \gamma \colon \R  \to \R $ and a continuously differentiable function $y \colon \R  \to \R $ that satisfy the following for all $ t \in \R $:

% $$\beta (t) = \int^1_0 D_1 f (sx(t),x(\eta_x(t))) \,\intd s,\qquad  \gamma (t) = \int^1_0 D_2 f (0,sx(\eta_x(t))) \,\intd s < 0$$
% \begin{equation}
% \label{felbontas}
% \dot{x}(t) = \beta(t)x(t) + \gamma(t)x(\eta_x(t))
% \end{equation}
% $$d(t) = \exp \left(- \int^t_{t_0} \beta(s) \,\intd s\right) > 0,\qquad  c(t) = \gamma(t) \exp \left(- \int^t_{\eta_x(t)} \beta(s) \,\intd s\right) < 0$$
% \begin{equation}
% \label{y-felbontas}
% y(t) = d(t) x(t),\qquad \dot{y}(t) = c(t) y( \eta_x (t) ), \qquad t \in \R , \end{equation}
% where $t_0\in\R $ is fixed arbitrarily.	
% \end{corollary}

% From this we have $\sgn y(t) = \sgn x(t)$. Letting  
% \begin{equation}\label{L1def}
% L_1\coloneqq \max_{\xi\in[-M,M]^2} \|Df(\xi)\|,
% \end{equation}
% we obtain, that inequalities
% \begin{align*}
% |\beta (t)| &\leq L_1,\qquad
% - L_1 \leq \gamma(t) < 0 ,\qquad
% e^{ - L_1 |t-t_0|} \leq d(t) \leq e^{L_1 |t-t_0|} \quad \mbox{and}\quad
% %\shortintertext{and}
% - L_1 e^{L_1 K}\leq c(t) < 0
% \end{align*}
% are fulfilled for all $t\in \R $.

\section{The discrete Lyapunov function}
\label{sec:signchange}

Let $\varphi \in C([a, b],\R  )\setminus \{ \0 \}$. 
%and $I \supset [a,b]$. 
We define the sign change counting function as follows:
\[\scc(\varphi , [a,b] ) = 0, \quad \text{if}\quad \varphi (s)\geq 0\quad \text{or}\quad \varphi (s)\leq 0\quad \text{for all } s \in  [a,b],
\]
and otherwise
\begin{equation}
\label{sc}
\begin{aligned} \operatorname{sc}(\varphi , [a,b])= \sup \{& k \in \mathbb {N}: k \ge 1, \text{ there } \text{ exists } (s^i)_0^k \subset [a,b], \text { such that} \\& s^0< s^1< \dots< s^k \text {, } \varphi (s^{i-1}) \varphi (s^i) < 0 \text { for all } 1 \le i \le k \}. \end{aligned}
\end{equation}
Now let
\begin{equation}
\label{tilde}
V(\varphi,[a,b]) =
\begin{cases}
\scc(\varphi,[a,b]), &\text{if it is odd or } \infty \text{,}\\
\scc(\varphi,[a,b]) + 1, &\text{otherwise.}
\end{cases}
\end{equation}
Obviously, $V(\varphi,[a,b]) \in \{1,3,\ldots\} \cup \{ \infty \}$.

%Furthermore, if $z \colon \Dom z \to \R $ is continuous and $\Dom z \subset \R $ is a closed interval, then let 
%\begin{equation}
%\label{sima-v}
%V(z) := V(z,\Dom z) \text{.}
%\end{equation}

%Furthermore, if $x$ is a solution of equation \eqref{eq:dde}, we will use notation
%\begin{equation}
%\label{sima-v}
%V(x_t) \coloneqq V(x,[\eta_x(t),t]).
%\end{equation}
%We adopt this notation also for solutions of equations \eqref{felbontas} and \eqref{y-felbontas}.

We define a class of regular functions on the interval $[a,b]$ in the following way:
\begin{equation}
%\begin{split}
%H_{[a,b]}= \bigl\{ \varphi \in C^1 ( [a,b] , \R  ) : {}& \varphi (b) \neq 0 \text{, or } \varphi(a)\dot{\varphi}(b) < 0 \text{,}\\
%&\varphi(a) \neq 0 \text{, or } \dot{\varphi}(a) \varphi(b) > 0 \text{,}\\
%&\text{$\varphi$ has only simple zeros} \bigr\}.
%\end{split}
\begin{aligned} \begin{aligned} H_{[a,b]}= \bigl \{& \varphi \in C^1 ( [a,b], \mathbb {R}): \varphi (b) \ne 0 \text {, or } \varphi (a)\dot{\varphi }(b) < 0 \text {,}\\&\varphi (a) \ne 0 \text {, or } \dot{\varphi }(a) \varphi (b) > 0 \text {, } \varphi \text { has only simple zeros} \bigr \}. \end{aligned} \end{aligned}
\end{equation}

In the following, we recall from \cite{krisztin-arino} some important properties of $V$. First, we describe the continuity of $V$ and its boundedness over regular functions. 

\begin{lemma}[{\cite[Lemma~4.1]{krisztin-arino}}]\leavevmode
\label{alapveto-v}
\begin{enumerate}[label=(\roman*)]
\item $V$ is semi-continuous from below in the following sense: if $\varphi$ and $\varphi^n$ are non-zero continuous functions on $[a,b]$ and $[a^n,b^n]$ for all $n\in\naturals$, respectively, and
$$\max_{s \in [a,b] \cap [a^n,b^n]} |\varphi(s)-\varphi^n(s)| \to 0,\qquad
a^n \to a\quad \mbox{and}\quad b^n \to b \quad \text{as } n \to \infty \text{,}$$
then
$$V(\varphi,[a,b]) \leq \liminf_{n \to \infty} V(\varphi^n,[a^n,b^n]) \text{.}$$
\item If $\varphi \in H_{[a,b]}$, then $V(\varphi,[a,b]) < \infty$.
\item If there exists $\delta > 0$ such that $\varphi \in C^1([a-\delta, b+\delta], \R )$ and
$\varphi |_{[a,b]} \in H_{[a,b]}$, then there exists $\gamma \in (0, \delta)$, such that if
$$|a-c| < \gamma,\qquad |b-d| < \gamma,\qquad \psi \in 
C^1([c,d], \R ), \quad \mbox{and} \quad \|\psi - \varphi\|_{C^1([c,d], \R )} < \gamma $$
hold, then
\[
V(\psi,[c,d]) = V(\varphi,[a,b]).
\]
\end{enumerate}
\end{lemma}

\medskip

%%%%%%%%%%%%%%%%%%%%%%%%%%%%%%%%%%%%%%%%%%%%%%%%%%%%%%%%%%%%%%%%%%%%%%%%
%\subsection*{Monotonicity of $\bm{\tilde V}$}

We consider equations of the form 
\begin{equation}\label{transformed_maineq}
\dot{y}(t) = c(t)y(t-\tau(t)),
\end{equation} 
akin to the auxiliary equation of \Cref{subsec:trafo}, 
and study the monotonicity of $V$ along delay-long solution segments. 
To that end, we introduce some notation and assumptions that will be applied throughout the remaining part of the section.

Let $I = [a,b]$, $c \colon I \to \R , \tau : I \to \R $ be continuous functions such that $c(t) < 0$, $\tau(t) > 0$ for all $t \in I$, and $\eta \colon I\to \R $, $\eta(t)= t - \tau(t)$ be strictly increasing on $I$. 

Set $k \in \N \setminus \{ 0, 1 \}$. Suppose that there is a finite sequence $(\alpha_j)^k_{j=1}$ in $[a,b]$ such that 
$\eta(\alpha_j) = \alpha_{j-1}$ if $j \in \{ 2, \ldots, k \}$. \added{Recall that $\eta^n$ denotes the $n$-th iterate of $\eta$.} 

Finally, set $J\coloneqq \{ t - \tau(t) : t \in I \} \cup I$ and let $y\colon J \to \R $ be a continuous function that is continuously differentiable on $I$ and satisfies \eqref{transformed_maineq} there.

\begin{lemma}[{\cite[Lemma~4.2]{krisztin-arino}}]
\label{monoton-v}
Suppose that $I = [a,b], c, \tau, \eta, y, k$ and $(\alpha_j)^k_{j=1}$ are as above, and $\left.y\right|_{[\eta(t),t]}$ is not zero for $t \in I$. Then, the following statements hold:
\begin{enumerate}[label=(\roman*)]
\item if $t^1, t^2 \in I$, $t^1 < t^2$, then $V(y, [\eta(t^2),t^2]) \leq V(y, [\eta(t^1),t^1])$,
\item if $k \geq 3$, $t \in [\alpha_3,b]$, and $y(t) = y(\eta(t)) = 0$, then $V(y, [\eta(t),t]) < V(y, [\eta^3(t),\eta^2(t)])$ or $V(y, [\eta(t),t]) = \infty$,
\item if $k \geq 4$, $t \in [\alpha_4,b]$, and $V(y, [\eta(t),t]) = V(y, [\eta^4(t),\eta^3(t)]) < \infty$, then $y|_{[\eta(t),t]} \in H_{[\eta(t),t]}$.
\end{enumerate}
\end{lemma}

The following lemma is a slight generalization of \cite[Lemma~4.3]{krisztin-arino} giving control on the decay of solutions of \eqref{transformed_maineq}. As we will see, the proof is analogous to the proof in \cite{krisztin-arino}.

\begin{lemma}
\label{korlatozo-v}
Fix an integer $n \geq 1$ and assume that the real numbers $t'<t$ are such that $c \colon [t',t] \to \R $ and $\tau \colon [t',t] \to \R $ are continuous, $\eta \colon [t',t]\ni s \mapsto s - \tau(s)$ is strictly increasing and $t'\leq \eta^{2n+4}(t)$.  Moreover, let $y$ be a continuous function on  $[\eta^{2n+5}(t), t]$ that is differentiable on $[\eta^{2n+4}(t), t]$ and satisfies \eqref{transformed_maineq} there.

Suppose also that there exist positive constants $c_0, c_1, \tau_0$ and $L_{\tau}$, such that
\begin{equation*}
\begin{split}
- c_1 & \leq c(s) \leq -c_0, \\
\tau_0 & \leq \tau(s), \\
|\tau(s^1) - \tau(s^2)| & \leq L_{\tau} |s^1 - s^2|
\end{split}
\end{equation*}
hold for all $s,s^1,s^2 \in [t',t]$. Then the following statements hold.
\begin{enumerate}[label=(\roman*)]
\item For any $\delta\in (0,\tau_0)$ \added{and positive integer $m\leq 2n+4$ there exists $C_m(\delta)$, independent of $t$ and $y$, such that for any closed interval $\Delta\subset [\eta^{m}(t),t]$ with length $\delta$
$$ \min_{s\in\Delta}|y(s)|\leq C_m(\delta)\max_{s\in [\eta(t),t]} |y(s)|$$
holds. }
\item If in addition $n$ is odd and  $V(y,[\eta^{2n+4}(t),\eta^{2n+3}(t)]) = n$, 
then there exists a constant $k = k(n, c_0, c_1, \tau_0, L_{\tau}) > 0$, independent of $t$ and $y$, such that
$$\max_{s \in [\eta^2(t),\eta(t)]} |y(s)| \leq k \max_{s \in [\eta(t),t]} |y(s)|$$
holds.
\end{enumerate}
\end{lemma}

\begin{proof}
The proof is similar to that of \cite[Lemma~4.3]{krisztin-arino} which deals with the special case $n=1$. For the readers' convenience, we include a detailed proof.

\medskip
\noindent
\added{In the following, let $\bar{y}=\max _{s \in[\eta(t), t]}|y(s)|$.}
\medskip

\noindent\added{\textbf{Part (i).} 
The statement is clearly valid for $m=1$ and arbitrary $\delta\in (0,\tau_0)$ with $C_1(\delta)=1$.}

\added{Now, let us consider $m=2$, and assume that $\Delta \subset [\eta^2(t), \eta(t)]$, that is, $\Delta=[\eta(s^1), \eta(s^2)]$, $\delta=\eta(s^2)-\eta(s^1)$ and $\eta(t) \leq s^1<s^2 \leq t$. Integrating \eqref{transformed_maineq} on $[s^1, s^2]$, we obtain}
\begin{equation}\label{eq:ys2-ys1}
\added{y(s^2)-y(s^1)=\int_{s^1}^{s^2} c(u) y(\eta(u)) \,d u}
\end{equation}
\added{The length of $[s^1, s^2]$ can be estimated by}
\begin{equation}\label{eq:delta-est}
\added{\delta=\eta(s^2)-\eta(s^1) \leq s^2-s^1+|\tau(s^2)-\tau(s^1)| \leq (1+L_\tau)(s^2-s^1).}
\end{equation}
\added{Combining equation \eqref{eq:ys2-ys1} with inequality \eqref{eq:delta-est} we gain
$$
\min _{s \in \Delta}|y(s)| \leq \frac{2(1+L_\tau)}{c_0 \delta} \bar{y}.
$$
This implies that the statement holds for $m=2$ with 
$C_2(\delta ):=\max \bigl \{1,\frac{2(1+L_\tau )}{c_0 \delta }\bigr \}$.}

\added{We will prove the statement for arbitrary $m\leq 2n+4$ by mathematical induction.  Assume now that the claim holds for some $m\leq 2n+3$. Consider any interval $\Delta \subset[\eta^{m+1}(t), t]$ of length $\delta$. If the length of $\Delta \cap[\eta^m(t), t]$ is greater than or equal to $\delta / 2$, then we choose $C_{m+1}(\delta)=C_m(\delta / 2)$. Assume that $|\Delta \cap[\eta^{m+1}(t), \eta^m(t)]|>\delta / 2$. There are $t^1, t^2 \in$ $[\eta^m(t), \eta^{m-1}(t)]$ such that $[\eta(t^1), \eta(t^2)] \subset \Delta$ and $\eta(t^2)-\eta(t^1)=\delta / 2$. From the Lipschitz continuity of $\tau$, we obtain
$$
t^2-t^1 \geqslant \frac{\delta}{2(1+L_\tau)}.
$$
Considering the intervals
$$
\left[t^1, t^1+\frac{\delta}{6\left(1+L_\tau\right)}\right],\left[t^2-\frac{\delta}{6\left(1+L_\tau\right)}, t^2\right] \subset\left[\eta^m(t), \eta^{m-1}(t)\right]
$$
of length $\bar{\delta}=\delta /(6\left(1+L_\tau\right))$, the first part of the proof gives that
$$
\min _{s \in \left[t^1, t^1+\bar{\delta}\right]}|y(s)| \leq C_m(\bar{\delta}) \bar{y}, \quad \text{and}\quad \min_{s \in\left[t^2-\bar{\delta}, t^2\right]}|y(s)| \leq  C_m(\bar{\delta}) \bar{y}.
$$
Applying the mean value theorem, we obtain a $t^* \in\left(t^1, t^2\right)$ such that
$$
\left|\dot{y}\left(t^*\right)\right| \leq \frac{2 C_m(\bar{\delta}) \bar{y}}{\bar{\delta}}.
$$
Using equation \eqref{transformed_maineq} we infer
$$
|y(\eta(t^*))| \leq \frac{|\dot{y}(t^*)|}{c_0} \leq \frac{2 C_m(\bar{\delta}) \bar{y}}{c_0 \bar{\delta}}.
$$
Since $\eta(t^*) \in \Delta$, it follows that
$$
\min _{s \in \Delta}|y(s)| \leq \frac{2 C_m(\bar{\delta})}{c_0 \bar{\delta}} \bar{y}.
$$
Then the statement holds for $m+1$  with 
$$
C_{m+1}(\delta)=\max \left\{C_m({\delta}/{2}), \frac{2C_m(\bar{\delta})}{c_0 \bar{\delta}}\right\}.
$$
This completes the proof of statement (i).}
\medskip

\noindent \textbf{Part (ii).} Let us fix a positive, odd $n$ and choose $\delta>0$ small enough that $2\delta(1+L_\tau)(2+L_\tau)^n<\tau_0$ holds and let \added{$C=C_{2n+4}(\delta)$ be given by statement (i)}. We may assume without loss of generality that $C>1$. Moreover, for integers $m$ with $0\leq m \leq n+1$ let $k_m\coloneqq C(1+\delta c_1)^{n+1-m}$. We claim that statement (ii) holds with $k\coloneqq k_0=C(1+\delta c_1)^{n+1}$.

% \begin{figure}\centering
%     \subcaptionbox{Illustration of the first step ($i=1$) in part (ii).\label{techlemma-fig2}}
% 	{\includegraphics[width=0.49\textwidth]{technical-lemma-fig2.pdf}}
%     \hfill
%     \subcaptionbox{Illustration of step $i=2$ in part (ii).\label{techlemma-fig1}}
% 	{\includegraphics[width=0.49\textwidth]{technical-lemma-fig1.pdf}}
	
% \caption{\added{Illustration of the proof of \cref{korlatozo-v}}.}
% 	\label{fig:techlemma}
% \end{figure}

\begin{figure}[htpb!]\centering
	{\includegraphics[width=0.8\textwidth]{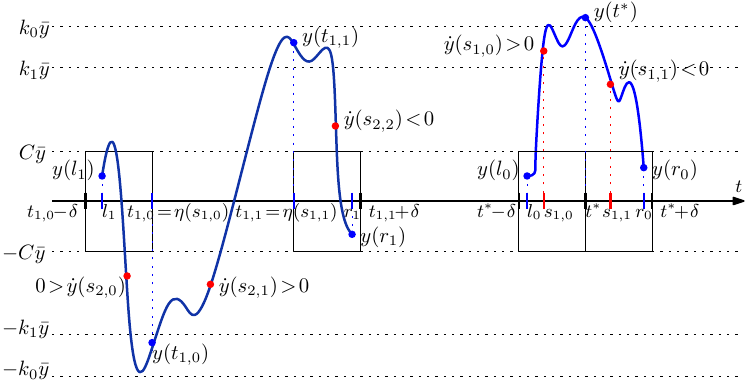}}
\caption{\added{Illustration of the proof of \cref{korlatozo-v}. On the right-hand side, we see how we can find two values ($s_{1,0}$ and $s_{1,1}$), using part (i) and the mean value theorem, where the derivatives have opposite signs. This results in a sign change of $y$ between $t_{1,0}=\eta(s_{1,0})$ and $t_{1,1}=\eta(s_{1,1})$. The latter is illustrated on the left-hand side of the figure. Then the procedure is continued to obtain 3 values ($s_{2,0}$, $s_{2,1}$ and $s_{2,2}$) with interchanging signs of derivatives of $y$.}	\label{fig:techlemma}}
\end{figure}
Let us indirectly assume that there exists $t^*\in [\eta^2(t),\eta(t)]$, such that $|y(t^*)|>k \bar{y}$. \added{The  proof goes as follows: starting from the interval $J\coloneqq [t^*-\delta,t^*+\delta]$ and using part (i) we get two points in the interval whose derivatives have opposite signs. In view of \eqref{transformed_maineq} we obtain that $y$ changes sign at least once in the interval $\eta(J)$, i.e.\ $\scc(y,\eta(J))\geq 1$. Then a similar argument is repeated altogether $n+2$ times to arrive at $\scc(y,\eta^{n+2}(J))\geq n+2$. The choice of $\delta$ guarantees that the length of the interval $\eta^{n+2}(J)$ is less than the minimal length of the delay ($\tau_0)$. This contradicts the assumption $V(y,[\eta^{2n+4}(t),\eta^{2n+3}(t)]) = n$. The first and second steps are illustrated in \cref{fig:techlemma}.}

\added{Now let us see the details.} For definiteness, suppose that $y(t^*)$ is positive (the negative case is analogous). Then by statement (i) we obtain that there exists $l_0\in [t^*-\delta,t^*]\subset [\eta^3(t),\eta(t)]$ and $r_0\in [t^*,t^*+\delta]\subset [\eta^2(t),t]$, such that $|y(l_0)|\leq C \bar{y}$ and $|y(r_0)|\leq C \bar{y}$ hold. Then one can apply the mean value theorem to deduce that there exist $s_{1,0}\in [l_0,t^*]\subset [\eta^3(t),\eta(t)]$ and $s_{1,1}\in [t^*,r_0]\subset [\eta^2(t),t]$, such that 
$$\dot{y}(s_{1,0}) > \frac{(k_0-C)\bar{y}}{\delta} \quad \text{and} 
\quad \dot{y}(s_{1,1}) < -\frac{(k_0-C)\bar{y}}{\delta}$$
hold. \added{This is illustrated on the right-hand side of  \cref{fig:techlemma}}. Then \eqref{transformed_maineq} implies that for $t_{1,j}\coloneqq \eta(s_{1,j})$ ($j=0,1$),
$$-y(t_{1,0})>\frac{(k_0-C)\bar{y}}{\delta c_1}\geq \frac{(k_0-k_1)\bar{y}}{\delta c_1}= k_1\bar{y}\quad \text{and}\quad y(t_{1,1})> \frac{(k_0-C)\bar{y}}{\delta c_1}\geq \frac{(k_0-k_1)\bar{y}}{\delta c_1}= k_1\bar{y} $$
are satisfied. Clearly, $t_{1,0}\in [\eta(t^*-\delta),\eta(t^*)]\subset[\eta^4(t),\eta^2(t)]$, while $t_{1,1}\in[\eta(t^*),\eta(t^*+\delta)]\subset [\eta^3(t),\eta(t)]$, and consequently $y$ has at least one sign change on $[\eta(t^*-\delta),\eta(t^*+\delta)]\subset [\eta^4(t),\eta(t)]$.

Now, let us define $a_0\coloneqq b_0\coloneqq t^*$ and for all $1\leq m\leq n+1$, $a_{m}\coloneqq \eta(a_{m-1}-\delta)$ and $b_m\coloneqq \eta(b_{m-1}+\delta)$, respectively. First of all, the monotonicity of $\eta$, and $\delta<\tau_0$ imply that $a_i$ and $b_i$ are well-defined, and that $[a_i,b_i]\subset [\eta^{2i+2}(t),\eta(t)]$ holds for all $0\leq i\leq n+1$.

We claim that our assumptions imply that for all $1\leq i\leq n+1$, there exists $t_{i,0}<t_{i,1}<\dots<t_{i,i}$, all from $[a_i, b_i]$, such that 
$(-1)^{i-j}y(t_{i,j})>0$ for all $0\leq j\leq i$, and moreover, $y(t_{i,i})>k_i\bar{y}$ and  $(-1)^i y(t_{i,0})>k_i\bar{y}$ also hold. 

We have just shown that these are fulfilled for $i=1$ (with same notations). \added{We will show that the claim holds for $i+1$, too. The step from $i=1$ to $i=2$ is illustrated on the left-hand side of  \cref{fig:techlemma}.}

Now assume that for some $1\leq i\leq n$, the claim holds. Then, by statement (i), there exist $l_i\in [t_{i,0}-\delta,t_{i,0}]$ and $r_i\in [t_{i,i},t_{i,i}+\delta]$, such that $|y(l_i)|\leq C\bar{y}$ and $|y(r_i)|\leq C\bar{y}$ hold.  Applying the mean value theorem  yields that there exist $s_{i+1,0}\in [l_i,t_{i,0}]$ and $s_{i+1,i+1}\in [t_{i,i},r_i]$, such that 
$$(-1)^i\dot{y}(s_{i+1,0}) > \frac{(k_i-C)\bar{y}}{\delta} \quad \text{and} 
\quad \dot{y}(s_{i+1,i+1}) < -\frac{(k_i-C)\bar{y}}{\delta}.$$
Then \eqref{transformed_maineq} implies that for $t_{i+1,0}\coloneqq \eta(s_{i+1,0})$ and $t_{i+1,i+1}\coloneqq \eta(s_{i+1,i+1})$,
\begin{align}
(-1)^{i+1}y(t_{i+1,0})&>\frac{(k_i-C)\bar{y}}{\delta c_1}\geq \frac{(k_i-k_{i+1})\bar{y}}{\delta c_1}= k_{i+1}\bar{y}\label{y_bigl}\\
\shortintertext{and}
y(t_{i+1,i+1})&> \frac{(k_i-C)\bar{y}}{\delta c_1}\geq \frac{(k_i-k_{i+1})\bar{y}}{\delta c_1} =k_{i+1}\bar{y}   \label{y_bigr}
\end{align}
are fulfilled. 

On the other hand, by application of the mean value theorem, we obtain that for all $1\leq j\leq i$ there exist $s_{i+1,j}\in (t_{i,j-1},t_{i,j})$ such that $(-1)^{i-j}\dot{y}(s_{i+1,j})>0$. Now, if we define $t_{i+1,j}\coloneqq \eta(s_{i+1,j})$ for all $1\leq j\leq i$, then by \eqref{transformed_maineq} we infer that $(-1)^{i+1-j}y(t_{i+1,j})>0$ for all $1\leq j\leq i$, which together with \eqref{y_bigl} and \eqref{y_bigr} and with the monotonicity of $\eta$ shows that the claim holds for $i+1$, too. This proves our claim.

It follows from the above claim  that $V(y,[a_{n+1},b_{n+1}])\geq n+2$. For a contradiction, we need  to prove that $b_{n+1}-a_{n+1}\leq \tau_0$. First notice that 
\begin{equation}\label{estimate-eta-diff}
0\leq \eta(t_2)-\eta(t_1)\leq (t_2-t_1)(1+L_\tau)
\end{equation}
holds for arbitrary $t'\leq t_1\leq t_2\leq t$. 

We state that $b_i-a_i\leq 2\delta(1+L_\tau)(2+L_\tau)^{i-1}$ is fulfilled for all $1\leq i\leq n+1$. The estimate \eqref{estimate-eta-diff} readily implies that
\[
b_1-a_1=\eta(t^*+\delta)-\eta(t^*-\delta)\leq 2\delta(1+L_\tau),
\]
hence the claim holds for $i=1$. Assume now that it holds for some $1\leq i\leq n$. Then from inequality \eqref{estimate-eta-diff} and from the inductive hypothesis it follows that
\begin{align*}
b_{i+1}-a_{i+1}&= \eta(b_i+\delta)-\eta(a_i-\delta)\\
&\leq  2\delta\bigl[(1+L_\tau)(2+L_\tau)^{i-1}+1\bigr](1+L_\tau)\\
&\leq 2\delta(1+L_\tau)(2+L_\tau)^i
\end{align*}
holds, which proves the statement. Consequently, $b_{n+1}-a_{n+1}<\tau_0$ holds.

This combined with $a_{n+1}\geq \eta^{2n+4}(t)$, $V(y,[a_{n+1},b_{n+1}])\geq n+2$ and with Lemma~\ref{monoton-v}\,(i) means that $V(y,[\eta^{2n+4}(t),\eta^{2n+3}(t))\geq n+2$, a contradiction to $V(y,[\eta^{2n+4}(t),\eta^{2n+3}(t])= n$. This concludes the proof of statement (ii).
% See also \cite[Section~2]{arino} for constant delay. See also Theorem 7.4 and Corollary 7.5 of \cite{mallet-paret} constant delay, only proved on the attractor.
\end{proof}

%\begin{remark}\label{rem:properties:of:V}
	Consider a nontrivial entire solution $x$ of \eqref{eq:dde} and let $\tau(t) \coloneqq r(x_t)$, $\eta(t) \coloneqq \eta_x(t)=t-\tau(t)$ and $y(t), c(t)$ be the functions constructed in \cref{felb}. As  $\sgn x(t) = \sgn y(t)$, therefore, the function $V$ attains the same value on any fixed interval for both $x$ and $y$, and thus the statements from \Cref{monoton-v,korlatozo-v} apply for $x$, as well.
	
	% In particular, \Cref{monoton-v}\,(i) yields that $t\mapsto V(x,[\eta_x(t),t]))$ is decreasing in $t$, that is, we have a ``discrete Lyapunov function'' at hand. In (iii) we need some simple calculations in addition to see that the appropriate segment of $x$ satisfies the regularity conditions as well.
% \end{remark}

\begin{lemma}
\label{v-properties-on-sol}
Suppose that $x$ is a nontrivial entire solution of \eqref{eq:dde} with standing assumptions {\rm\ref{eq:f:bounded}--\ref{eta:increasing}} and let $\eta\coloneqq \eta_x$. Then the following statements hold:
\begin{enumerate}[label=(\roman*)]
\item if $t^1, t^2 \in I, t^1 < t^2$, then $V(x, [\eta(t^1),t^1]) \geq V(x, [\eta(t^2),t^2])$,
\item if  $x(t) = x(\eta(t)) = 0$ for some $t\in \R $, then $V(x, [\eta(t),t]) = \infty$ or $V(x, [\eta(t),t]) < V(x, [\eta^3(t),\eta^2(t)])$,
\item if $V(x, [\eta(t),t]) = V(x, [\eta^4(t),\eta^3(t)]) < \infty$ for some $t\in \R $, then $x|_{[\eta(t),t]} \in H_{[\eta(t),t]}$.
\item For all odd integer $n\in\N$, there exists $k>0$  such that if $x$ is a nontrivial entire solution with $V(x,[\eta^{2n+4}(t),\eta^{2n+3}(t)]) = n$, 
then 
$$\max_{s \in [\eta^2(t),\eta(t)]} |x(s)| \leq k \max_{s \in [\eta(t),t]} |x(s)|$$
holds for all $t\in \R$.
\end{enumerate}
\end{lemma}
\begin{proof} Let $y$ be defined by \eqref{eq:x_to_y_transform}. Then $y$ is an entire solution of \eqref{y-felbontas}, and \cref{monoton-v} can be applied to $y$. %It is also clear that $\sgn y(t)=\sgn x(t)$ for all $t$.

Statements (i) and (ii) are straightforward consequences of \cref{monoton-v}\,(i)--(ii) and \cref{felb}. 

Statement (iv) follows directly from \cref{felb} and \cref{korlatozo-v}\,(ii).

The proof of statement (iii) follows the argument found in \cite[pp.~212--213]{krisztin-walther-wu}, with minor changes. For the reader's convenience, we give the details here. We need to show that
	$$x|_{[\eta_x(t),t]}\in H_{[\eta_x(t),t]}\quad\text{holds, whenever}\quad y_{[\eta_x(t),t]}\in H_{[\eta_x(t),t]}.$$
	To see this, let $t\in \R$ such that $y_{[\eta_x(t),t]}\in H_{[\eta_x(t),t]}$ is satisfied. Then $y$ is continuously differentiable and
	\begin{equation}
		\dot y(s)=\exp\left(-\int_{0}^sa(\theta)d\theta\right)\left(\dot x(s)-a(s)x(s)\right)\label{eq:doty}
	\end{equation}
	for $s\in \R$.
 
 Assume that $s\in(\eta_x(t),t)$ and $x(s)=0$. From equation \eqref{eq:x_to_y_transform} and $y|_{[\eta_x(t),t]}\in H_{[\eta_x(t),t]}$ it follows that $y(s)=0$ and $\dot y(s)\ne0$. Hence, \eqref{eq:doty} implies $\dot x(s)\ne 0$ and $x|_{[\eta_x(t),t]}\in H_{[\eta_x(t),t]}$. 
 
 Now assume that $x(t)=0$. Then $y_{[\eta_x(t),t]}\in H_{[\eta_x(t),t]}$ and \eqref{eq:x_to_y_transform} imply $y(t)=0$ and $\dot y(t)y(\eta_x(t))<0$. Hence, \eqref{eq:doty} with  $s=t$ yields $\dot x(t) \dot y(t)>0$. Equations \eqref{felbontas} and \eqref{y-felbontas} combined with $x(t)=0$ and the negativity of functions $b$ and $c$  imply that $x(\eta_x(t))y(\eta_x(t))>0$. Thus $\dot x(t)x(\eta_x(t))<0$ and $x|_{[\eta_x(t),t]}\in H_{[\eta_x(t),t]}$. 
 
 Finally, assume that  $x(\eta_x(t))=0$. It follows from $y|_{[\eta_x(t),t]}\in H_{[\eta_x(t),t]}$ and \eqref{eq:x_to_y_transform} that $y(\eta_x(t))=0$ and $\dot{y}(\eta_x(t))y(t)>0$. This combined with \eqref{eq:x_to_y_transform} imply $x(t)y(t)>0$. By virtue of equation \eqref{eq:doty} with $s=\eta_x(t)$ we obtain $\dot{x}(\eta_x(t))\dot{y}(\eta_x(t))>0$. Now, combining these three inequalities, we obtain $\dot x(\eta_x(t))x(t)>0$, yielding $x|_{[\eta_x(t),t]}\in H_{[\eta_x(t),t]}$.

To prove statement (iv), fix an arbitrary $t_0\in \R$, and define $y$ and $c$  by transformation \eqref{eq:x_to_y_transform}. Applying \cref{korlatozo-v}\,(ii) -- and using its notation -- we infer 
\begin{align*}\max_{s\in [\eta^2(t_0),\eta(t_0)]} |x(s)| 
&\leq e^{2Ka_0}\max_{s\in [\eta^2(t_0),\eta(t_0)]} |y(s)| \\
&\leq ke^{2Ka_0} \max_{s\in [\eta(t_0),t_0]} |y(s)|\\
&\leq ke^{3Ka_0} \max_{s\in [\eta(t_0),t_0]} |x(s)|,
\end{align*}
as required.
\end{proof}
% Thanks to hypotheses  \ref{eq:delay:bounded}, \ref{eq:f:bounded} and \ref{eq:delay:lip}, we may define
% \[\tau_0\coloneqq \min_{\varphi\in \overline \phasespace} r(\varphi)\qquad L_1\coloneqq \max_{\xi\in [-M,M]^2}|D_1 f(\xi)|+\max_{\xi\in [-M,M]^2}|D_2f(\xi)|. \]
%%%%%%%%%%%%%%%%%%%%%%%%%%%%%%%%%%%%%%%%%%%%%%%%%%%%%%%%%%%%%%%%%%%%%%%%
\section{The Morse decomposition}
\label{sec:morse}

This section is devoted to the proof of our main result, \cref{thm:morse}. In order to state the theorem, first we need to introduce some notation and the definition of a Morse decomposition itself.

\begin{definition}
\label{def:m1m2}
	A \textit{Morse decomposition} of the global attractor $\A$ is a finite, ordered system $\M_0 \prec \M_1 \prec \dots \prec \M_m$ of nonempty, compact, invariant, pairwise disjoint subsets of $\A$, which is ordered in the following sense: 
for all $\varphi \in \A$ and any  entire solution $x$ with $x_0=\varphi$, there exist $0\leq j\leq i\leq m$ with
\begin{enumerate}[label=$(m_{\arabic*})$, leftmargin=*]
	\item $\alpha(x)\subseteq \M_i$ and $\omega(\varphi)\subseteq \M_j$, \label{morse_prop1}
	\item $i=j$ implies $\varphi \in \M_i$ (thus, $x_t\in \M_i$ for all $t\in\R$). \label{morse_prop2}
\end{enumerate}

Although the Morse sets should be -- by definition -- nonempty, for the simplicity of notation we take the liberty to occasionally list a finite number of sets $\M_0,\dots,\M_m$ as a Morse decomposition with some (but not all) of them being potentially empty. This should not cause any confusion.
\end{definition}

Note that -- similarly as in \cite{mallet-paret} and \cite{polner} -- due to the lack of backwards uniqueness, we need to deal with $\alpha$-limit sets of concrete entire solutions, rather than that of an element of $\phasespace$. This is reflected in our notation.

In order to state the main result of this section, we need to study the following linear equation with constant delay, associated to the flow  $\semiflow$ and to equation \eqref{eq:dde} (see \cite{cooke-huang}):
\begin{equation}\label{linearized_eq}
\dot{x}(t)=A x(t) + B x(t-\tau),
\end{equation}
where $A= D_1 f(0,0)$, $B=D_2 f(0,0)<0$ and $\tau=r(\0)> 0$. The eigenvalues of \eqref{linearized_eq} are exactly those $\lambda\in\complex$ numbers that solve the characteristic equation 
\begin{equation}\label{char_eq}
\lambda-A - e^{-\lambda\tau}B=0.
\end{equation}

Let us denote by $M^\ast$  the number of eigenvalues (counting multiplicities) of  \eqref{linearized_eq} with strictly positive real part. This is known to be finite \cite[Theorem~6.1]{mallet-paret}. Moreover, let
\begin{align*}
%N^\ast_{+}&\coloneqq\begin{cases}
%M^\ast+1,&\mbox{if the origin is nonhyperbolic and $M^\ast$ is odd},\\
%M^\ast,&\mbox{otherwise}. 
%\end{cases}\\
N^\ast &\coloneqq\begin{cases}
M^\ast+1,&\mbox{if the origin is nonhyperbolic and $M^\ast$ is even},\\
M^\ast,&\mbox{otherwise}. 
\end{cases}
\end{align*}

It is worth mentioning that provided $ A + B < 0$, then $M^\ast$ is always even  by \cite[Theorem~6.1]{mallet-paret}.

Let us define for $N \in \N\setminus\{ N^*\}$ the sets
\begin{equation}\label{morse:sets:normal}
\begin{aligned} \begin{aligned} \mathcal {S}_N :=\Bigl \{\varphi \in \mathcal {A}\setminus \{ \textbf{0} \}: &\text{ there } \text{ exists } \text{ an } \text{ entire } \text{ solution }  x  \text{ with } x_0=\varphi \text{, } \\&\text{ such } \text{ that } V(x,[\eta _x(t),t]) = N \text{ for } \text{ all } t \in \mathbb {R} \text{ and } \textbf{0} \notin \alpha (x) \cup \omega (\varphi ) \Bigr \},\\ \mathcal {S}_{N}^+ :=\Bigl \{\varphi \in \mathcal {A}\setminus \{ \textbf{0} \}: &\text{ there } \text{ exists } \text{ an } \text{ entire } \text{ solution }  x  \text{ with } x_0=\varphi \text{, } \\&\text{ such } \text{ that } V(x,[\eta _x(t),t]) \ge N \text{ for } \text{ all } t \in \mathbb {R} \text{ and } \textbf{0} \notin \alpha (x) \cup \omega (\varphi ) \Bigr \}. \end{aligned} \end{aligned}
%\S_N \coloneqq \Bigl\{\varphi \in \Anull : {}&\mbox{there exists an entire solution $x$ with $x_0=\varphi$,}	\\
%& \mbox{such that } V(x,[\eta_x(t),t]) = N \mbox{ for all } t \in \R  \mbox{ and } \0 \notin \alpha(x) \cup \omega(\varphi) \Bigr\},\\
%\S_{N}^+ \coloneqq \Bigl\{\varphi \in \Anull : {}&\mbox{there exists an entire solution $x$ with $x_0=\varphi$,}	\\
%& \mbox{such that } V(x,[\eta_x(t),t]) \geq N \mbox{ for all } t \in \R  \mbox{ and } \0 \notin \alpha(x) \cup \omega(\varphi) \Bigr\}.
%\end{aligned}
\end{equation}
Furthermore, let
%\begin{center}
%	$\S_{N^*} := \{ \0 \}$, if the origin is hyperbolic
%\end{center}
%otherwise
\begin{equation}\label{morse:sets:N*}
\begin{aligned} \begin{aligned} \mathcal {S}_{N^*} :=\Bigl \{\varphi \in \mathcal {A}\setminus \{ \textbf{0} \}: &\text{ there } \text{ exists } \text{ an } \text{ entire } \text{ solution }  x  \text{ with } x_0=\varphi \text{, } \\&\text{ such } \text{ that } V(x,[\eta _x(t),t]) = N^* \text{ for } \text{ all } t \in \mathbb {R}\Bigr \} \cup \{ \textbf{0} \}. \end{aligned} \end{aligned}
%\begin{aligned}
%\S_{N^*} \coloneqq \Bigl\{\varphi \in \Anull : {}&\mbox{there exists an entire solution $x$ with $x_0=\varphi$,}	\\
%& \mbox{such that } V(x,[\eta_x(t),t]) = N^\ast \mbox{ for all } t \in \R \Bigr\}  \cup \{ \0 \}.
%\end{aligned}
\end{equation}
%In any other case let $\S_N$ undefined.

Note that if $ A + B < 0$ and the origin is hyperbolic, then $\S_{N^\ast}=\{\0\}$.

\bigskip
Now we are in position to state our main result.

\begin{theorem}
	\label{thm:morse}
	Suppose that hypotheses {\rm\ref{eq:f:bounded}--\ref{eta:increasing}} hold and $N_0>N^\ast$ is an odd number. For $n\leq N_0$ let
 \[\M_n=\begin{cases}
     \S_n, &\text{if }  n<N_0 ,\\
     \S_{N_0}^+, &\text{if } n=N_0.
 \end{cases}\]
Then the sets $\M_0,\dots,\M_{N_0}$ form a Morse decomposition of the global attractor of \eqref{eq:dde}.
\end{theorem}

%The proof of \Cref{thm_morse} follows the argument seen in \cite{polner}. First one shows that there is an upper bound $\Vmax\in \N$, such that   $V(x,[\eta_x(t),t])\leq \Vmax$ for all $t\in\R $ and any entire solution $x$ of \eqref{eq:dde} -- see \cref{v-korl}. Although the main ideas of its proof is analogous to that of \cite{polner}, we provide a detailed proof of it. The reason for that, apart from striving to be self-contained,  is that the proofs reveal some technical difficulties compelling us to restrict our study to those two types of delays introduced in \Cref{subsec:delays}. %they contain a proof of the important property: bounded growth backwards in time.

% Having \cref{v-korl} at hand, it is left to show that the Morse sets $\S_N$ are compact (invariance and disjointness are clear from definition), and that the Morse property \ref{morse:property} is satisfied. The proof can be carried out via a number of lemmas, analogous to Lemmas~4.3--4.8 in \cite{polner}, with only minor changes, therefore we shall merely state those and omit the details here. We mention that a more detailed argument showing how those lemmas imply the theorem is presented in \cite{garab-poetzsche} and in \cite{garab} for analogous results on difference equations and on cyclic system of differential equations with constant delay.

The proof follows the logic of \cite{polner}, however, there are several differences. On the one hand, there are technicalities arising due to the state-dependency of the delay. On the other hand, with this generality of the delay, we are not able to prove boundedness of the Lyapunov function on the global attractor. The sets of type $\S_{N_0}^+$ are introduced to overcome this difficulty. These sets contain all ``rapidly oscillating'' entire solutions (i.e.~along which $V\geq N_0$ holds) and potentially super-high frequency solutions (i.e.~that change sign infinitely many times on a finite interval), should they exist. Nevertheless, in \cref{sec:examples} we consider two major types of state-dependent delays that allow us to prove the boundedness of the discrete Lyapunov function and to omit sets of rapid oscillation.

We need several lemmas to prove the theorem. The first one utilizes the fact that the eigenvalues are nicely ordered (see \cite[Theorem 6.1]{mallet-paret}). It also reveals the role of $N^*$.

\begin{samepage}\begin{lemma}\leavevmode\label{lemma_boundedsols_lin}
	\begin{enumerate}[label=(\roman*)]
		\item If $x$ is a nontrivial bounded solution of \eqref{linearized_eq} on $(-\infty,0]$, then $V(x,[t-\tau,t]))\leq N^\ast$ for all  $t\leq 0$.
		\item If $x$ is a nontrivial bounded solution of \eqref{linearized_eq} on $[0,\infty)$, then $V(x,[t-\tau,t]))\geq N^\ast$ for all  $t\geq 0$.	
	\end{enumerate}
\end{lemma}
\end{samepage}
\begin{proof}
This is a special case of \cite[Propositions~3.8]{garab}.
\end{proof}

The next lemma is a key technical result stating, roughly speaking, that the above property of the linearized equation is preserved locally in the nonlinear equation.
%  build a bridge between our main equation \eqref{eq:dde} and its linearization \eqref{linearized_eq}, and

\begin{lemma}
\label{u-korny}
There exists  an open neighborhood $\U$ of $\0$ in $\A$ such that if
$x \colon \R  \to (-M,M)$ is a nontrivial entire solution of \eqref{eq:dde}, then
\begin{enumerate}[label=(\roman*)]
\item if $x_t \in \overline{\U}$ for all  $t \leq 0$, then 	$V(x,[\eta_x(t),t]) \leq N^\ast$ holds for all $t \in \R $.
\item if $x_t \in \overline{\U}$, for all $t \geq 0$, then $V(x,[\eta_x(t),t]) \geq N^\ast$ holds for all $t \in \R $.
\end{enumerate}
\end{lemma}

\begin{proof}
(i) We present an indirect argument: suppose that there exists $(\varphi^k)^\infty_{k=0} \subset \Anull$ so that $ \sup_{t\leq 0}\|x^k_{t}\| \to 0$, as $k \to \infty$, and $V(x^k,[\eta_k(t_k),t_k])>N^*$ holds for all $k\in \N_0$ for some $t_k\in \R $, where $x^k \coloneqq x^{\varphi^k}$ and $\eta_k\coloneqq \eta_{x^k}$.

Let $s_k\leq 0$, such that $|x^k(s_k)| > \frac{1}{2} \sup_{t \leq 0}\|x^k_{t}\|$. By virtue of \cref{v-properties-on-sol}\,(i) and because \eqref{eq:dde} is autonomous, we may assume that $|x^k(0)| > \frac{1}{2} \sup_{t \leq 0}\|x^k_{t}\|$ and $V(x^k, [\eta_k(t),t])> N^*$ holds for all $t \leq 0$.

According to \cref{felb}, $x^k$ solves the  nonautonomous linear  equation
\[\dot{x}^k(t) = a^k(t) x^k(t) + b^k(t) x^k(\eta_{k}(t)),\qquad t \in \R ,\]
where $a^k$ and $b^k$ are defined by
\[a^k (t) = \int^1_0 D_1 f (sx^k(t),x^k(\eta_k(t))) \,\intd s,\qquad  
b^k (t) = \int^1_0 D_2 f (0,sx^k(\eta_k(t))) \,\intd s < 0.\]
Since $\|x^k_t\| \to 0$ for all $t \leq 0$ as $k\to\infty$, $f$ is continuously differentiable, and $r\colon \phasespace\to \R $ is Lipschitz continuous,  
\begin{equation}\label{conv:to:ABtau}
\lim_{k\to\infty}a^{k}(t)= A=D_1f(0,0),\qquad  \lim_{k\to\infty}b^{k}(t)= B=D_2f(0,0), \quad \mbox{and}\quad  \lim_{k\to\infty}\eta_k(t)= t-\tau
\end{equation}
hold for all $t\leq 0$, where $\tau\coloneqq r(\0)$.

Let
\[z^k(t) \coloneqq \frac{x^k(t)}{\|x^k_0\|},\qquad t \in \R .\]
Then
\begin{equation}
\label{z-bontas}
\dot{z}^k(t) = a^k(t) z^k(t) + b^k(t) z^k(\eta_{k}(t)),\qquad  t \in \R 
\end{equation}
is satisfied as well, and we know that
\[V(z^k,[\eta_{k}(t),t]) > N^\ast,\]
and $\|z^k_t\| < 2$, if $t \leq 0$ and $\|z^k_0\| = 1$. Thus $(z^k)_{k=0}^\infty$ is uniformly bounded, and equation \eqref{z-bontas} together the bounds $a_0$ and $b_1$ from \cref{felb} imply that it is also equicontinuous. So we can apply the Arzel\`a--Ascoli theorem and the Cantor diagonalization to obtain that there is a subsequence $(z^{k_l})_{l=0}^\infty$ and a continuous function $z\colon (-\infty,0]\to \R$ such that
\begin{equation}\label{z:conv}
z^{k_l}(t) \to z(t),\quad \mbox{as } l\to \infty
\end{equation}
holds for all $t\leq 0$, and convergence is uniform in $t$ on any compact subset of $(-\infty,0]$. It can also be easily seen that the convergence is uniform  in \eqref{conv:to:ABtau} on compact subsets of $(-\infty,0]$. These together with equations \eqref{z-bontas} and \eqref{z:conv} imply that $z$ is differentiable, $\dot{z}^{k_l}\to \dot{z}$, as $l\to \infty$, and convergence is uniform on compact subsets of $(-\infty,0]$, moreover, $z$ satisfies equation \eqref{linearized_eq}
%\[\dot{z}(t)=A z(t) + Bz(t-\tau)\]
on the interval $(-\infty,0]$. 

From the definition of functions $z^k$ it follows that $\|z_t\|\leq 2$ for all $t\leq 0$, and that $\|z_0\| = 1$. In particular, $z$ is a nontrivial solution of \eqref{linearized_eq} that is bounded on $(-\infty, 0]$, thus \cref{lemma_boundedsols_lin} yields that $V(z,[t-\tau,t]))\leq N^\ast$ holds for all  $t\leq 0$.
 
Note that  there exists $T < 0$, such that $z|_{[T-\tau,T]} \in H_{[T-\tau,T]}$. This can be seen by applying  transformation \eqref{eq:x_to_y_transform} to \eqref{linearized_eq} and then \cref{monoton-v}\,(iii). Finally, using the local uniform $C^1$-convergence of $(z^{k_l})_{l=0}^\infty$ and 	\cref{alapveto-v}\,(iii) we obtain that
\[N^\ast < % \lim_{l \to \infty} V(x^{k_l},[\eta_{k_l}(T),T]) = 
\lim_{l \to \infty} V(z^{k_l},[\eta_{k_l}(T),T]) = V(z,[T-\tau,T]) \leq N^\ast.\]
This is a contradiction.
\medskip

(ii) Arguing again by contradiction, let us suppose that there exists a sequence $(\varphi^k)^\infty_{n=1} \subset \Anull$ and entire solutions $x^k$ through $\varphi^k$ such that 
\begin{equation}\label{sup_xk_goes_to_0}
\sup_{t \in [0,\infty)}\|x^k_{t}\| \to 0,\quad \mbox{as } k \to \infty,
\end{equation}
and $V(x^k, [\eta_k(t_k),t_k]) < N^*$ is satisfied for some $t_k \in \R $, where $\eta_k \coloneqq \eta_{x^k}$. 

Let $m_0\in \N$ such that $m_0\tau_0\geq 2K$, where $\tau_0 = \min_{\varphi \in \A}r(\varphi)$, which is positive by hypotheses \ref{eq:delay:bounded}--\ref{eq:delay:lip}. We aim to apply \cref{v-properties-on-sol}\,(iv) ($m_0-1$ times) to arrive at a contradiction. Keeping this in mind, put  $m\coloneqq 2n+4+m_0$. 

% and let $\delta_0 > 0$ be so small that
% \[| D_2f(0,v) - B| \leq \frac{|B|}{2}, \quad \mbox{for all } |v| \leq \delta_0\]
% holds. 

Using the autonomy of equation \eqref{eq:dde}, the monotonicity of $V$ (see \cref{v-properties-on-sol}), and  assumption \eqref{sup_xk_goes_to_0} we may assume (by also passing to a subsequence if necessary) that there  exist positive integers $n$ and $k_0$ such that for all $k\geq k_0$  we have 
\begin{alignat}{3}
V(x^k,[\eta(t),t])&=n < N^*, &&\quad \mbox{if } t \geq -mK, \\
\shortintertext{and}
\label{bounded_from_above_by_2}
\|x^k_{0}\| &> \frac{1}{2} |x^k(t)|, &&\quad \mbox{if }  t\geq -K.	
% &\|x^k_{t}\| \leq \min \left\{ \frac{\tau}{2 \Lip{r}} , \delta_0 \right\},&&\quad  \mbox{if }  t\geq -mK,	
\end{alignat}

Since $\varphi^k \neq \0$, hence $x^k_t \neq \0$ for all $t \in \R $.  Applying \cref{v-properties-on-sol}\,(iv)  $m_0-1$ times, we infer the existence of a constant $C>1$ such that inequality
\begin{equation*}
\max_{s\in [\eta^{m_0}_k(t_0),\eta_k(t_0)]} |x^k(s)| 
\leq C^{m_0-1} \max_{s\in [\eta_k(t_0),t_0]} |x^k(s)|
\end{equation*}
holds for all $k\geq k_0$ and $t_0\geq 0$. Combining this with the definition of $m_0$, we obtain that also
\begin{equation}\label{exp_estimate_xk}
\max_{s\in [t_0-2K,\eta_k(t_0)]} |x^k(s)| \leq C^{m_0-1} \max_{s\in [\eta_k(t_0),t_0]} |x^k(s)|
\end{equation}
holds for all $k\geq k_0$ and $t_0\geq 0$.

Now we will proceed similarly as in part (i) and introduce the notations $z^k, a^k$ and $b^k$.
% . Let us define
% \[z^k(t) \coloneqq \frac{x^k(t)}{\|x^k_0\|}, \qquad t \in \R .\]
Putting together estimates \eqref{sup_xk_goes_to_0} and \eqref{exp_estimate_xk}   we obtain that the estimate
\[
|z^k(t)|\leq \max\bigl\{2,C^{m_0-1}\bigr\}
\]
holds for all $t\geq -2K$ and $k\geq k_0$. 

Now we can deduce the same way as in part (i) that there exist locally uniformly convergent subsequences
\[z^{k_l} \to z,\quad \mbox{and}\quad \dot{z}^{k_l} \to \dot{z}.\]
Moreover, since 
\[\dot{z}^{k_l}(t) = a^{k_l}(t) z^{k_l}(t) + a^{k_l}(t) z^{k_l}(\eta_{k_l}(t))\]
holds for all $t\geq 0$, it follows that the function $z$ satisfies the linear equation \eqref{linearized_eq} on $[0,\infty)$. Note also that since $\|z^{k_l}(0)\|=1$ for all $k\geq k_0$, hence $\|z_0\| = 1$ holds. Therefore, $z$ is a nontrivial solution of equation \eqref{linearized_eq} that is bounded on $[0,\infty)$, so we can apply \cref{lemma_boundedsols_lin} to obtain that $V(z,[t-\tau,t]) \geq N^\ast$ for all $t\in [0,\infty)$, and thanks to the monotonicity of $V$ this also holds for all $t \in \R $.

Finally, applying transformation \eqref{eq:x_to_y_transform} to solution $z$ and then \cref{monoton-v}\,(iii), we obtain the existence of a $T> 0$, such that $z|_{[T-\tau,T]} \in H_{[T-\tau,T]}$. Then, using the local uniform $C^1$-convergence of $(z^{k_l})_{l=0}^\infty$ and 	\cref{alapveto-v}\,(iii) we arrive at 
\[N^\ast > \lim_{l \to \infty} V(x^{k_l},[\eta_{k_l}(T),T]) = \lim_{l \to \infty} V(z^{k_l},[\eta_{k_l}(T),T]) = V(z,[T-\tau,T]) \geq N^\ast,\]
which is a contradiction.
\end{proof}

The next three lemmas are crucial in proving both the compactness of the Morse sets and that the Morse properties \ref{morse_prop1}--\ref{morse_prop2} \added{(as in  Definition~\ref{def:m1m2})} hold.

\begin{lemma}\label{C1-norm-lemma}
    Suppose that $(\varphi^n)_{n=0}^\infty\subset \A$ and $\varphi\in \A$ with \added{$\lim_{n\to \infty}\|\varphi^n - \varphi\|=0$}. Then there exists a subsequence  $(\varphi^{n_k})_{n=0}^\infty$ that converges to $\varphi$ in the $C^1$-norm, i.e. $\lim_{k\to \infty} \|\varphi^{n_k} - \varphi\|_1=0$.
\end{lemma}
\begin{proof}
The proof of \cite[Lemma~4.3]{polner} for the constant delay case (with positive feedback) applies with trivial modifications.
\end{proof}

\begin{samepage}
\begin{lemma}\label{lemma:V:on:limit:sets}
\leavevmode
\begin{enumerate}[label=(\roman*)]
    \item Let $\varphi \in \Anull$,  $N\in \N$  such that $\lim_{t\to \infty} V(x^{\varphi}_t, [\eta_{x^{\varphi}}(t),t])=N$, and $y$ be an entire solution through an arbitrary $\psi\in \omega(\varphi)\setminus \{\0\}$. Then $V(y,[\eta_y(0),0])=N$.
    \item Let $\varphi \in \Anull$,  $N\in \N\cup \{\infty\}$  such that $\lim_{t\to \infty} V(x^{\varphi}_t, [\eta_{x^{\varphi}}(t),t])\geq N$, and $y$ be an entire solution through an arbitrary $\psi\in \omega(\varphi)\setminus \{\0\}$. Then $V(y,[\eta_y(0),0]\geq N$.
%    \item Let $\varphi \in \Anull$,  $\lim_{n\to \infty} V(x^{\varphi}_t, [\eta_{x^{\varphi}}(t),t])= \infty$, and $y$ an entire solution through an arbitrary $\psi\in \omega(\varphi)\setminus \{\0\}$. Then $V(y,[\eta_y(0),0])= \infty$.
    \item Let $\varphi \in \Anull$,  $N\in \N$ such that $\lim_{t\to -\infty} V(x^{\varphi}_t, [\eta_{x^{\varphi}}(t),t])=N$, and $y$ be an entire solution through an arbitrary $\psi\in \alpha(x)\setminus \{\0\}$. Then $V(y,[\eta_y(0),0])=N$.
    \item Let $\varphi \in \Anull$,  $N\in \N\cup \{\infty\}$  such that $\lim_{t\to -\infty} V(x^{\varphi}_t, [\eta_{x^{\varphi}}(t),t])\geq N$, and $y$ be an entire solution through an arbitrary $\psi\in \alpha(x)\setminus \{\0\}$. Then $V(y,[\eta_y(0),0]\geq N$.
\end{enumerate}
\end{lemma}
\end{samepage}
\begin{proof}
%    The proof is similar to that of \cite[Lemma~4.4]{polner}.

Let $\psi\in \omega(\varphi)\setminus\{0\}$. In view of \cref{C1-norm-lemma}  there exists $t_n\to \infty$ such that 
\begin{equation}\label{xtn_to_psi}
\lim_{n\to \infty} \|x_{t_n}-\psi\|_1 = 0.
\end{equation}
By invariance, for any entire solution $y$ through $\psi\in \omega(\varphi)\setminus\{0\}$, $y_t\in \omega(\varphi)$ holds for all $t\in \R $.
Thus, for any $t\in \R $, there exists $t'_n\to \infty$ with $\lim_{n\to \infty} x_{t'_n} = y_{t}$. 

Note that this combined with \cref{v-properties-on-sol}\,(i) implies that $V(y,[\eta_y(t),t]) \leq N$ holds for all $t\in\R $ in case the assumptions of statement (i) are imposed.

Thanks to \cref{v-properties-on-sol}\,(i), there exists $K\coloneqq \lim_{t\to \infty} V(y,[\eta_y(t),t])$. We claim that $K\geq N$.

To see this, let us indirectly assume $K<N$. This yields in particular that $K$ is finite. Then, it follows from \cref{v-properties-on-sol}\,(iii) that there exists $T_2>0$ such that $V(y, [\eta_y(T_2),T_2])=K$ and $y|_{[\eta_y(T_2),T_2]}\in H_{[\eta_y(T_2), T_2]}$. Note that  $\lim_{n\to \infty}\|x_{t_n+T_2}-y_{T_2}\|_1 = 0$ holds, 
thanks to \eqref{xtn_to_psi}. Now, applying \cref{alapveto-v}\,(iii) and our assumption on $x$, we obtain
\begin{equation*}
N\leq \lim_{n\to \infty} V(x,[\eta_x(t_n+T_2), t_n+T_2])=V(y,[\eta_y(T_2),T_2])=K,
\end{equation*}
a contradiction. 

Hence, $V(y,[\eta_y(0),0])\geq K\geq N$ holds by virtue of \cref{v-properties-on-sol}\,(i).

This readily proves statement (i), as $V(y,[\eta_y(0),0]) \leq N$ has already been shown above.

As for statement (ii), it follows from $K\geq N$ and the monotonicity of $V$ that  $V(y,[\eta_y(0),0]\geq N$ holds, as stated.

Statements (iii) and (iv) can be proved analogously.
% (i) This combined with \cref{v-properties-on-sol}\,(i) implies that $V(y,[\eta_y(t),t]) \leq N$. Furthermore, from \cref{v-properties-on-sol}\,(iii) follows that there exists $T_1<0<T_2$ such that $y|_{[\eta_y(T_i),T_i]}\in H_{[\eta_y(T_i), T_i]}\ (i=1,2)$. In view of \cref{C1-norm-lemma}, there exists $t_n'\to \infty$ such that $\lim_{n\to \infty}\|x_{t_n'}-y_{T_1}\|_1 = 0$. On the other hand, thanks to \eqref{xtn_to_psi}, $\lim_{n\to \infty}\|x_{t_n+T_2}-y_{T_2}\|_1=0$, so, by applying \cref{alapveto-v}\,(iii), we obtain
% \begin{align*}V(y,[\eta_y(T_1),T_1]&=\lim_{n\to \infty} V(x,[\eta_x(t_n'),t_n'])\\
% &=N=\lim_{n\to \infty} V(x,[\eta_x(t_n+T_2), t_n+T_2])=V(y,[\eta_y(T_2),T_2]).\end{align*}
% Using the monotonicity of $V$ (i.e.\ \cref{v-properties-on-sol}\,(i)) we arrive at $V(y,[\eta_y(t),t])=N$ for all $t\in [T_1,T_2]$, so, in particular, $V(y,[\eta_y(0),0])=N$, as stated.
\end{proof}

\begin{lemma}\label{lemma:0:is:equal:or:not-in:limit:set}
Let $\varphi \in \Anull$ and $x$ be an entire solution with $x_0=\varphi$. Then the following statements hold.
\begin{enumerate}[label=(\roman*)]
\item If  $ \lim_{t \to \infty}V(x,[\eta_x(t),t])\neq N^\ast$, then either $\omega(\varphi)=\{\0\}$ or else $\0\notin \omega(\varphi)$.
\item If $ \lim_{t \to -\infty}V(x,[\eta_x(t),t])\neq N^\ast$, then either $\alpha(x)=\{\0\}$ or else $\0\notin \alpha(x)$.
\end{enumerate}
\end{lemma}

\begin{proof} The proof is analogous to that of \cite[Lemma 4.5]{polner}.

(i) Let us assume that $\0\in \omega(\varphi)\neq \{\0\}$. It suffices to prove that $\lim_{t \to \infty}V(x,[\eta_x(t),t])= N^\ast$. 

Let $\U$ denote the neighbourhood of $\0$ from \cref{u-korny}.  Then there exists $\U_1\subseteq \U$ such that $\0 \in \U_1$ and $\omega(\varphi)\not\subseteq \overline{\U_1}$. This means that $x$ enters and leaves $\U_1$ infinitely many times as $t\to \infty$.  More precisely, there exist positive sequences $(t_n), (\sigma_n)$ and $(\tau_n)$, such that $t_n \nearrow \infty$ and $\lim_{n\to\infty}x_{t_n}=\0$, moreover,  the intervals $(t_n-\sigma_n, t_n+\tau_n)$ are pairwise disjoint (i.e.\ $t_n+\tau_n < t_{n+1}-\sigma_{n+1}),\ \forall n\in\N)$, $x_t\in \U_1$ for all $t\in (t_n-\sigma_n, t_n+\tau_n)$ and $x_{t_n-\sigma_n},x_{t_n+\tau_n} \in \partial \U_1$ hold for all $n\in \N$. Note that from \ref{eq:delay:bounded} and \cref{v-properties-on-sol}\,(iv), it follows that $\lim_{n\to\infty}\sigma_n =\infty$.

Now, let 
\begin{align*}
y^n\colon \R  \ni t&\mapsto x(t_n-\sigma_n+t)\in \R ,\\
z^n\colon \R  \ni t&\mapsto x(t_n+\tau_n+t)\in \R .
\end{align*}
Similarly as in the proof of \cref{u-korny} we obtain -- by virtue of the Arzel\`a--Ascoli theorem and Cantor's diagonalization -- that there exist subsequences (still denoted by $y^n$ and $z^n$) and entire solutions $y$ and $z$ of \eqref{eq:dde} such that
\[y^n\to y,\quad \dot y^n \to \dot y\quad \text{and}\quad z^n\to z,\quad \dot z^n \to \dot z,\]
uniformly on compact subsets of $\R$, as $n \to \infty$. 

It follows that $y_t\in \overline{\U}_1$ for $t\geq 0,\ y_0\in \partial \U_1$ and $z_t\in \overline{\U}_1$ for $t\leq 0$. Using the monotonicity of $V$ (\cref{v-properties-on-sol}\,(i)) and \cref{u-korny} we infer
\[V(y,[\eta_y(t),t])\geq N^\ast\geq V(z, [\eta_z(t),t]),\quad \text{for all } t\in \R.\]
In particular, we have $V(y,[\eta_y(0),0])\geq N^\ast\geq V(z, [\eta_z(0),0])$. On the other hand, as $y_0=\lim_{n\to \infty}x_{t_n-\sigma_n}$, and similarly $z_0=\lim_{n\to \infty}x_{t_n+\tau_n}$ holds, therefore $y_0,z_0 \in \omega(\varphi)\setminus \{\0\}$, hence, in view of \cref{lemma:V:on:limit:sets},
\[\lim_{t\to \infty} V(x,[\eta_x(t),t]) =V(y,[\eta_y(0),0]) =  V(z, [\eta_z(0),0])
\]
holds. This implies $\lim_{t\to \infty} V(x,[\eta_x(t),t])= N^\ast$, which proves our claim.

The proof of statement (ii) is analogous.
\end{proof}

\begin{lemma}\label{lemma:limit:set:trivial:or:in:Morse-set}
Let $N$ be a positive odd number, $\varphi \in \Anull$ and $x$ be an entire solution with $x_0=\varphi$, and let $\eta=\eta_x$. Then the following statements hold.
\begin{enumerate}[label=(\roman*)]
\item If $\lim_{t\to \infty} V(x, [\eta(t),t])=N\neq N^\ast$, then either $\omega(\varphi)=\{\0\}$ or $\omega(\varphi)\subseteq \S_N$.
\item If $\lim_{t\to \infty} V(x, [\eta(t),t])\geq N>N^\ast$, then either $\omega(\varphi)=\{\0\}$ or $\omega(\varphi)\subseteq \S_N^+$.
\item If $\lim_{t\to -\infty} V(x, [\eta(t),t])=N\neq N^\ast$, then either $\alpha(x)=\{\0\}$ or $\alpha(x)\subseteq \S_N$.
\item If $\lim_{t\to -\infty} V(x, [\eta(t),t])\geq N>N^\ast$, then either $\alpha(x)=\{\0\}$ or $\alpha(x)\subseteq \S_N^+$.
\end{enumerate}
\end{lemma}

\begin{proof}
Using the previous lemmas, the proof \cite[Lemma~4.6]{polner} applies almost verbatim for the state-dependent case, too. 

Let us prove statements (i) (and (ii)) simultaneously. Assume that $\omega(\varphi)\neq \{ \0 \}$ and \linebreak $\lim_{t\to \infty} V(x, [\eta(t),t])=N \neq N^\ast$ ($\lim_{t\to \infty} V(x, [\eta(t),t])\geq N > N^\ast$). Then $\0 \notin \omega(\varphi)$ follows by \cref{lemma:0:is:equal:or:not-in:limit:set}.

Let $\psi\in \omega(\varphi)$ and $y$ be an entire solution through $\psi$. Then, by invariance and compactness of omega-limit sets, $y_t \in \omega(\varphi)$ holds for all $t\in \R$, moreover, 
\[\alpha(y)\cup \omega(\psi) \subseteq \omega(\varphi),\]
and thus $\0 \notin \alpha(y)\cup \omega(\psi)$. 

Now, it follows from \cref{lemma:V:on:limit:sets} that $V(y,[\eta_y(0),0])=N$ ($V(y,[\eta_y(0),0])\geq N$) for all $t\in \R$, and hence, $\psi \in \S_N$ ($\psi \in \S_N^+$) holds.
\end{proof}

The next lemma shows that all Morse sets but $S_{N^*}$ are bounded away from $\0$.

\begin{lemma}\label{lemma:v-korny}
Let $N$ be a positive odd number. Then the following statements hold.
\begin{enumerate}[label=(\roman*)]
\item  If $N\neq N^\ast$, then there exists an open neighborhood $\V$ of $\0$ in $\phasespace$ such that $\V  \cap \S_N=\emptyset$.
\item  If $N> N^\ast$, then there exists an open neighborhood $\V$ of $\0$ in $\phasespace$ such that $\V \cap \S_N^+=\emptyset$. 
\end{enumerate}
\end{lemma}

\begin{proof} Let us consider the cases $N>N^\ast$ and treat statements (i) (in this special case) and (ii) simultaneously.

Let us assume to the contrary that there exists a sequence $(\varphi^n)_{n=1}^\infty$ in $\S_N\setminus \{\0\}$ (in $\S_N^+ \setminus \{\0\}$) such that $\lim_{n\to\infty}\varphi^n=\0$. 

Using notation $x^n=x^{\varphi^n}$, we have $V(x^n,[\eta_{x^n}(t),t])=N$ ($V(x^n,[\eta_{x^n}(t),t])\geq N>N^\ast$) for all $t\in \R$, furthermore, $\0 \notin \alpha(x^n)\cup \omega(\varphi^n)$ holds by \cref{lemma:0:is:equal:or:not-in:limit:set}.

Let $\U$ be the set from \cref{u-korny} and suppose indirectly that there exists $n\in \N$ such that $\omega(\varphi^n) \in \overline\U$. Then for an entire solution $y$ within $\omega(\varphi^n)$, $V(y,[\eta_y(t),t]) = N^\ast$ must hold for all $t\in\R$, by mentioned lemma. On the other hand, \cref{lemma:V:on:limit:sets} implies $V(y,[\eta_y(0),0])=N\neq N^\ast$ ($V(y,[\eta_y(0),0])\geq N> N^\ast$). Hence, $\omega(\varphi^n)\not\subseteq \overline\U$ holds for all $n\in \N$.

This implies the existence of a sequence $\sigma_n \to \infty$, such that for any $n$ large enough, $x^n_t\in \U$ for $t\in [0,\sigma_n)$ and $x^n_{\sigma_n} \in \partial \U$.

Now, similarly as in previous proofs, one may find a subsequence of $y^n\colon t \mapsto x^n(t+\sigma_n)$ (without changing notation) and an entire solution $y$, such that $y^n \to y$ and $\dot y^n \to \dot y$, as $n\to \infty$, uniformly on compact subsets of $\R$.

It follows by construction that  $y_t \in \overline\U$ for $t\leq 0$ and $y_0\in \partial U$ (in particular, $y$ is a nontrivial solution). By virtue of \cref{u-korny} and \cref{v-properties-on-sol}\,(i), we infer that $V(y,[\eta_y(t),t])\leq N^\ast$ holds for all $t\in \R$. Furthermore, \cref{v-properties-on-sol}\,(iii) and \cref{C1-norm-lemma} guarantee the existence of a $T>0$ and a subsequence (still denoted by $y^n$), for which
\[y|_{[\eta_y(T),T]}\in H_{[\eta_y(T),T]} \quad \text{and}\quad 
\lim_{n\to\infty}\|y_T^n - y_T\|_1=0.
\]
In light of \cref{alapveto-v}\,(iii) and the construction of $y$, we obtain 
\begin{equation*}
    \begin{aligned}  N^*<N&=\lim _{n\rightarrow \infty } V(x^n,[\eta (\sigma _n+T),\sigma _n+T])\\ & = \lim _{n\rightarrow \infty }V(y^n,[\eta _{y^n}(T),T])=V(y,[\eta _y(T),T])\le N^*, \end{aligned}
\end{equation*}
%\[N^\ast<N=\lim_{n\to\infty} V(x^n,[\eta(\sigma_n+T),\sigma_n+T])=
%\lim_{n\to\infty}V(y^n,[\eta_{y^n}(T),T])=V(y,[\eta_y(T),T])\leq N^\ast,\]
a contradiction that proves the statements.

It remains to prove statement (i) in case of $N<N^\ast$. The argument is completely analogous to the above presented one, with the only essential difference that, instead of the $\omega$-limit sets of the sequence $\varphi^n$, one has to deal with corresponding $\alpha$-limit sets.
\end{proof}

Now we can prove that the Morse sets are compact -- being closed subsets of the compact attractor $\A$.

\begin{lemma}\label{lemma:morse:sets:closed}
The sets $\S_{N^\ast}$, $\S_N$ and $\S_{N_0}^+$ are closed for  any odd numbers $N$ and $N_0$ with $N_0>N^\ast$.
\end{lemma}

\begin{proof}
The proof for $\S_N$ and $\S_{N^\ast}$ can be carried out by straightforward minor modifications of the proof of \cite[Lemma~4.8]{polner}, therefore it is omitted here.

Let us assume that $N_0>N^\ast$ is an odd number and $\varphi^n\in \S_{N_0}^+$ for all $n\in \N$ with $\lim_{n\to\infty}\varphi^n =\varphi\in \A$, and let us denote $x^{\varphi^n}$ by $x^n$. Then $V(x^n,[\eta_{x^n}(0),0])\geq N_0$ holds for all $n\in \N$.

As we have seen before, by using the Arzel\`a--Ascoli theorem and Cantor's diagonalization, we may assume that there is an entire solution $x$ with $x_0=\varphi$ such that   
\[x^n \to x, \qquad \dot x^n \to \dot x\]
uniformly on compact subsets of $\R$. We must show that $ V(x,[\eta_x(t),t])\geq N_0$ for all $t\in\R$.

Assume to the contrary that there exists $t_0\in\R$ with $V(x,[\eta_x(t_0),t_0])<N_0$. Then \cref{v-properties-on-sol} guarantees the existence of $T>0$ for which 
\begin{equation*}
\begin{aligned} 
\lim_{t\rightarrow \infty } V(x,[\eta_x(t),t])=V(x,[\eta _x(T),T]) \quad \text {and} \quad x|_{[\eta_x(T),T]}\in H|_{[\eta_x(T),T]} 
\end{aligned}
\end{equation*}
%\[\lim_{t\to \infty} V(x,[\eta_x(t),t]=V(x,[\eta_x(T),T] \quad \text{and} \quad x|_{[\eta_x(T),T]}\in H|_{[\eta_x(T),T]}\]
holds. According to \cref{C1-norm-lemma}, we may assume that $\lim_{n\to\infty}\|x_T^n-x_T^{\vphantom n}\|_1=0$ and thus, by virtue of \cref{alapveto-v}\,(iii) we infer
\[
N_0>V(x,[\eta_x(T),T])=V(x^n,[\eta_{x^n}(T),T])\geq N_0.
\]

This contradiction proves our statement.
\end{proof}

After all this preparatory work, we can prove our main result.

\begin{proof}[Proof of \cref{thm:morse}]
The argument follows ideas of the proofs presented in \cite[Theorem B]{mallet-paret}, \cite[Theorem~4.1]{garab-poetzsche} and \cite[Theorem 4.1]{garab}. The main difference here is the presence of the set $\M_{N_0}=\S_{N_0}^+$.
	 
By definition, the sets $\M_n$ are finitely many, they are invariant, pairwise disjoint, and, according to \Cref{lemma:morse:sets:closed}, they are compact. Hence, it is only left to prove that the Morse properties \ref{morse_prop1}--\ref{morse_prop2} hold.%, that is, for all $\varphi \in \A$ and any bounded entire solution $x$ for which $x_0=\varphi$ holds, there exist $i\geq j$ with $\alpha(x)\subseteq \M_i$ and $\omega(\varphi)\subseteq \M_j$, and in case $i=j$, then $\varphi \in \M_i$ (thus, $x_t\in \M_i$ for every $t\in \R$).

%This can be shown by combining \Cref{V_monotone,u-korny,lemma:0:is:equal:or:not-in:limit:set,lemma:limit:set:trivial:or:in:Morse-set,lemma:V:on:limit:sets}. For a detailed argument we refer the reader to \cite[Theorem~3.2]{garab}, where the existence of an analogous Morse decomposition was proved for delay difference equations.

Note that thanks to the backward uniqueness of the zero solution, \ref{morse_prop1}--\ref{morse_prop2} hold trivially for $\varphi=0$. So consider an arbitrary $\varphi\in \A\setminus \{\0\}$, and let $x$ be an entire solution of \eqref{eq:dde} for which $x_0=\varphi$ holds. Furthermore, define
\[i\coloneqq \lim_{t\to -\infty}V(x,[\eta_x(t),t]) \quad \mbox{and}\quad j\coloneqq \lim_{t\to \infty}V(x,[\eta_x(t),t]) .\]
From \cref{v-properties-on-sol}, one obtains that $j\leq i$. Note that $i$ or even $j$ may be $\infty$.

First, observe that if $j=N^\ast$, then $\omega(\varphi)\subseteq \M_{N^\ast}$. To see this, choose an arbitrary $\psi\in \omega(\varphi)$ and let $y$ be an entire solution through $\psi$. If $\psi=0$, then $\psi\in \M_{N^\ast}$ holds by definition, so we may assume  that $\psi\neq 0$. By \Cref{lemma:V:on:limit:sets} we infer  $V(y,[\eta_y(0),0])=N^\ast$. Moreover, by the invariance of $\omega(\varphi)\setminus \{\0\}$ and by virtue of \Cref{lemma:V:on:limit:sets}, $V(y,[\eta_y(t),t])=N^\ast$ holds for all $t\in \R$, and, in particular, $\psi \in \M_{N^\ast}$. A similar argument can be applied to prove that $i=N^\ast$ implies that $\alpha(x)\subseteq \M_{N^\ast}$ holds.

We will distinguish four cases in terms of the values of $i$ and $j$.
\medskip

\noindent \textit{Case 1.} If $i=j=N^\ast$, then $\alpha(x)\cup \omega(\varphi) \subseteq \M_{N^\ast}$ holds by the above observation. Moreover, from the monotonicity of $V$ it follows that $V(x,[\eta_x(t),t])\equiv N^\ast$ on $\R$, thus $x_t\in \M_{N^\ast}$ for all $t\in\R$, and both \ref{morse_prop1} and \ref{morse_prop2} hold.%
\medskip 

\noindent\textit{Case 2.} If $i>j=N^\ast$, then $\omega(\varphi)\subseteq \M_{N^\ast}$ holds. On the other hand, $\alpha(x)\neq \{\0\}$ holds. Otherwise, \Cref{u-korny} would imply $V(x, [\eta_x(t),t])\leq N^\ast$ for $t\in \R$, and thus $i\leq N^\ast=j$, which is a contradiction. Hence, \Cref{lemma:limit:set:trivial:or:in:Morse-set} yields $\alpha(x)\subseteq \M_i$ in case $i< N_0$, and $\alpha(x)\subseteq \M_{N_0}=S_{N_0}^+$, when $i\geq N_0$. Thus, property \ref{morse_prop1} is fulfilled. Note that \ref{morse_prop2} holds automatically, as the two Morse sets in question, i.e.\ $\M_{N^\ast}$ and $\M_i$, are different.
\medskip

\noindent\textit{Case 3.} A similar argument applies in the case when $i=N^\ast>j$.

\medskip
\noindent\textit{Case 4.} In case $i\neq N^\ast \neq j$, let 
\[
i_0 = \begin{cases}
    i, &\text{if } i<N_0,\\
    N_0, &\text{otherwise,}
\end{cases}
\qquad 
j_0 = \begin{cases}
    j, &\text{if } j<N_0,\\
    N_0, &\text{otherwise.}
\end{cases}
\]
Then \Cref{lemma:limit:set:trivial:or:in:Morse-set} yields that either $\omega(\varphi)=\{\0\}$ or $\omega(\varphi)\subseteq \M_{j_0}$. Similarly, either $\alpha(x)=\{\0\}$ or $\alpha(x)\subseteq \M_{i_0}$ holds. Note that $\omega(\varphi)$ and $\alpha(x)$ cannot be both $\{\0\}$ in this case, because then \Cref{u-korny} would imply that $V(x,[\eta_x(t),t])= N^\ast$ for all  $t\in \R$, a contradiction.

If none of  $\omega(\varphi)$ and  $\alpha(x)$ is $\{\0\}$, then from \Cref{lemma:limit:set:trivial:or:in:Morse-set} we obtain that $\omega(\varphi)\subseteq \M_{j_0}$ and $\alpha(x)\subseteq \M_{i_0}$ hold, so \ref{morse_prop1} is fulfilled. Furthermore, if $i_0=j_0$, then the definition of $\M_{i_0}$ and $\M_{j_0}$ imply that $V(x,[\eta_x(t),t])=i=j$ for all $t \in \R$, in case $i_0=j_0<N_0$, and $V(x,[\eta_x(t),t])\geq N_0$, when $i_0=j_0=N_0$. On the other hand, \Cref{lemma:0:is:equal:or:not-in:limit:set} ensures that $\0\notin \alpha(x)\cup\omega(\varphi)$, thus $x_t\in \M_i$  holds for all $t\in \R$. This establishes property  \ref{morse_prop2}.

If $\omega(\varphi)=\{\0\}\neq \alpha(x)$, then $\omega(\varphi)\subseteq \M_{N^\ast}$ holds by definition. Furthermore, \Cref{u-korny} implies that $V(x, [\eta_x(t),t])\geq N^\ast$ holds for all $t\in \R$, and consequently $N^\ast< j\leq i$. On the other hand, \Cref{lemma:limit:set:trivial:or:in:Morse-set} yields that $\alpha(x)\subseteq \M_{i_0}$, so \ref{morse_prop1} holds. Property \ref{morse_prop2} is fulfilled automatically. 

An analogous argument applies for the case when $\omega(\varphi)\neq \{0\}=\alpha(x)$.
\medskip

We have listed all possible cases, so the proof is complete.
\end{proof}

%%%%%%%%%%%%%%%%%%%%%%%%%%%%%%%%%%%%%%%%%%%%%%%%%%%%%%%%%%%%%%%%%%%%%%%%

\section{Boundedness of the Lyapunov function} \label{sec:sharper:results}
In this section we show that in case solutions fulfill the technical assumption \ref{iterated_zeros_imply_fullzero} below, then the discrete Lyapunov function has an upper bound on the global attractor. This allows for the following sharper version of \cref{thm:morse}.

\begin{theorem}
	\label{thm:morse:bdd}
	Suppose that hypotheses {\rm\ref{eq:f:bounded}--\ref{eta:increasing}} hold, and furthermore, it holds that
 \begin{enumerate}[label=$(\mathrm{H}_\infty)$,  leftmargin=*]
     \item \label{iterated_zeros_imply_fullzero} for any entire solution $x \colon \R  \to (-M,M)$ of \eqref{eq:dde}, if there exists $\sigma \in [\eta_x(0),0]$ such that $x(\eta^k_x(\sigma)) = 0$ for all $k \in \N$, then $x(t) = 0$ for all $t\in \R $.
 \end{enumerate}
 Then there exists $\Vmax\in \N$ such that   $V(x,[\eta_x(t),t])\leq \Vmax$ for all $t\in\R $ and any entire solution $x$ of \eqref{eq:dde}. Moreover, the sets $\M_0,\M_1,\dots, M_{V_{\max}}$ form a Morse decomposition of the global attractor of \eqref{eq:dde} with  $\M_n=\S_n$, $n=0,1, \dots, V_{\max}$.
\end{theorem}

\begin{figure}[htpb!]
    \centering
    \includegraphics[width=0.75\linewidth]{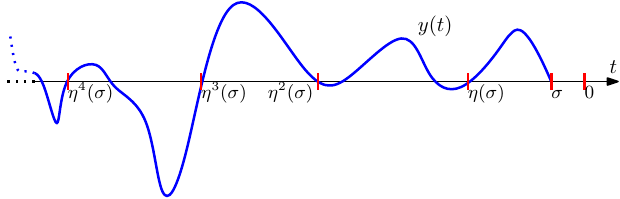}
    \caption{\added{Example of an entire solution not obeying condition \ref{iterated_zeros_imply_fullzero}.}}
    \label{fig:H-infty}
\end{figure}

\added{We do not know whether bounded entire solutions exists or not so that 
\ref{iterated_zeros_imply_fullzero} does not hold under conditions {\rm\ref{eq:f:bounded}--\ref{eta:increasing}}.
\cref{fig:H-infty} illustrates what scenarios are \textit{excluded} by hypothesis \ref{iterated_zeros_imply_fullzero}.}
In the next section, we will show two classes of state-dependent delays, for which our basic hypotheses {\rm\ref{eq:f:bounded}--\ref{eta:increasing}}, as well as \ref{iterated_zeros_imply_fullzero} are fulfilled and thus \cref{thm:morse:bdd} applies. 

To prove the above theorem, we need the following technical lemma. 

\begin{lemma}
\label{csak-nul}
Let hypotheses {\rm\ref{eq:f:bounded}--\ref{eta:increasing}}  and \ref{iterated_zeros_imply_fullzero} hold and $(\varphi^n)^\infty_{n=0} \subset \Anull$ be such that $\varphi^n \to \varphi \in \A$, and  $V(\varphi^n, [\eta_n(0),0]) \to \infty$,  as $n \to \infty$, where  $\eta_n\coloneqq \eta_{x^{\varphi^n}}$. Then $\varphi = \0$.
 % any  entire solution $x$ with $x_0=\varphi$, there exists $\sigma \in [\eta(0),0]$ such that $x(\eta^k(\sigma)) = 0$ for all $k \in \N$, where $\eta\coloneqq \eta_x$.
\end{lemma}

\begin{proof}  In view of \ref{iterated_zeros_imply_fullzero}, it clearly suffices to show that for the entire solution $x\coloneqq x^\varphi$, there exists $\sigma \in [\eta(0),0]$ such that $x(\eta^k(\sigma)) = 0$ for all $k \in \N$, where $\eta = \eta_x$. For brevity, let us introduce notation $x^n\coloneqq x^{\varphi^n}$. The proof is similar to that of \cite[Lemma~3.5]{polner}.
\medskip

\noindent\textit{Case 1:} $\varphi$ has finitely many zeros on the interval $[\eta(0),0]$: $\{\sigma_1, \dots,  \sigma_k\}$. Taking into account that $r(\varphi^n)\to r(\varphi)$ holds as $n\to \infty$, one can choose $\sigma \in \{\sigma_1, \dots, \sigma_k\}$ and a subsequence, still denoted by $(\varphi^n)^\infty_{k=0}$, in a way that 
\begin{equation}\label{sc_to_infty}
\scc(\varphi^{n},[\sigma-\varepsilon,\sigma+\varepsilon] \cap [-K,0]) \to \infty \quad \text{as } n \to \infty
\end{equation}
is satisfied for all $ \varepsilon > 0$.

Since $(x^{n})^\infty_{n=0}$ is uniformly bounded and, by \eqref{eq:dde} and \ref{eq:f:bounded}, so is $(\dot{x}^{n})^\infty_{n=0}$, therefore we can apply the Arzel\`a--Ascoli theorem followed by a Cantor diagonalization to obtain a subsequence, still denoted by $(x^{n})^\infty_{n=0}$, such that $x^{n}$ converges to a continuous function $\tilde{x}$, and convergence is uniform on any compact subset of $\R $. Bearing in mind that functions $x^{n}$ are solutions of \eqref{eq:dde}, we infer that $\tilde{x}$ is continuously differentiable and locally uniformly $\dot{x}^{n} \to \dot{\tilde{x}}$. % Theorem 7.17 of Rudin
Thus $x = \tilde{x}$ must be satisfied and $x(\sigma)=\varphi(\sigma)=0$ holds by our assumptions. 

In view of \eqref{sc_to_infty} and since $\dot{x}^{n}(\sigma) \to \dot{x}(\sigma)$, as $n\to \infty$, an indirect argument easily shows that $\dot{x}(\sigma) = 0$ also holds. Feedback condition \ref{eq:f:negative} yields that $x(\eta(\sigma)) = 0$. We claim that $\dot{x}(\eta(\sigma)) = 0$ is satisfied as well. We show this by an indirect argument:
suppose that there exists $\delta_1 > 0$ such that 
\begin{equation}\label{dotx>0}
%|\dot{x}(t)| > \frac{1}{2} |\dot{x}(\eta(\sigma))| > 0,\qquad t \in [\eta_x(\sigma) - \delta_1, \eta_x(\sigma) + \delta_1].
\begin{aligned} |\dot{x}(t)|> \frac{1}{2} |\dot{x}(\eta (\sigma ))| > 0,\qquad t \in [\eta (\sigma ) - \delta _1, \eta (\sigma ) + \delta _1]. \end{aligned}
\end{equation}
Then, by uniform convergence of $(x^{n})_{l=0}^\infty$, there exists $N_0$, that for $l \geq N_0$,
\begin{equation}\label{dotxn>0}
%|\dot{x}^{n}(t)| > \frac{1}{4} |\dot{x}(\eta(\sigma))| > 0,\qquad t \in [\eta_x(\sigma) - \delta_1, \eta_x(\sigma) + \delta_1].
\begin{aligned} |\dot{x}^{n}(t)|> \frac{1}{4} |\dot{x}(\eta (\sigma ))| > 0,\qquad t \in [\eta (\sigma ) - \delta _1, \eta (\sigma ) + \delta _1]. \end{aligned}
\end{equation}

Now, by the triangle inequality and from hypothesis \ref{eq:delay:lip} we obtain the estimates
\begin{align*}
|\eta_{n}(t)-\eta(\sigma)|&= |t-r_{x^{n}_t}-\sigma+r_{x_\sigma}|\\
&\leq |t-\sigma|+ \Lip{r}\|x^{n}_t - x_\sigma\|\\
&\leq |t-\sigma| + \Lip{r}(\|x^{n}_t - x_t\| + \|x_t-x_\sigma\|).
\end{align*}

Thus, local uniform convergence of $(x^n)_{n=0}^\infty$, hypothesis \ref{eq:f:bounded},  and monotonicity of the functions $\eta$ and $\eta_n$  imply that there exist
$N_1 \geq N_0$, $\delta_2 > 0$, such that if $l \geq N_1$, then:
\begin{align*}
[\eta(\sigma - \delta_2), \eta(\sigma + \delta_2)] &\subseteq [\eta(\sigma) - \delta_1, \eta(\sigma) + \delta_1],\\
\shortintertext{and}
[\eta_{n}(\sigma - \delta_2), \eta_{n}(\sigma + \delta_2)] &\subseteq  [\eta(\sigma) - \delta_1, \eta(\sigma) + \delta_1].
\end{align*}

Using these and inequality \eqref{dotxn>0}, it follows that for all $l\geq N_1$,  $x^{n	}$ has at most one zero in $[\eta_{n}(\sigma - \delta_2), \eta_{n}(\sigma + \delta_2)])$. According to \cref{felb}, the corresponding function $y^{n}$ has also maximum one zero there,  provided $l\geq N_1$. Therefore, $\dot{y}^{n}$ has at most one zero in the interval $[\sigma - \delta_2,\sigma + \delta_2]$, and consequently $y^{n}$ has at most two of that there, if $l\geq N_1$. Then, again by \cref{felb}, this holds for $x^{n}$ as well. This is a contradiction, thus  $\dot{x}(\eta(\sigma)) = 0$.

In a similar manner, one can deduct from \cref{felb} that for all $\varepsilon>0$,
\[\scc(x^{n},[\eta(\sigma) - \varepsilon , \eta(\sigma) + \varepsilon]) \to \infty, \quad \text{as }n \to \infty.\]
By iterating this argument 	we obtain that $x(\eta^k(\sigma)) = \dot{x}(\eta^k(\sigma)) = 0$ for all $k \in \N$.
\medskip

\noindent\textit{Case 2:} $\varphi$ has infinitely many zeros on $[\eta(0),0]$.
Let $\sigma$ be one of the accumulation points.
The function $y$ obtained from \cref{felb} shows that $\eta(\sigma)$ is an accumulation point of zeros of $y$ and $x$ too. This can be repeated for all $\eta^k(\sigma)$, $k \in \N$. Then continuity of $x$ yields that $x(\eta^k(\sigma))=0$ for all $k\in\N$.
\end{proof}

Now we can prove \cref{thm:morse:bdd}.

\begin{proof}[Proof of Theorem \ref{thm:morse:bdd}]
Having \cref{thm:morse} at hand, it is sufficient to prove the boundedness of the Lyapunov function on the global attractor. 

Suppose indirectly that there exists $(\varphi^n)^\infty_{n=0}$ and entire solutions $x^n$ of \eqref{eq:dde} such that $V(x^n, [\eta_n(0),0]) \to \infty$, as $n\to \infty$, where, as usual, $x^n\coloneqq x^{\varphi^n}$ and $\eta_n\coloneqq \eta_{x^n}$.

Since	 $\A$ is compact, we may assume that $\varphi^n$ converges to some $\varphi \in \A$, as $n\to\infty$. Then \cref{csak-nul} yields $\varphi = \0$.

This, together with $V(x^n, [\eta_n(0),0]) \to \infty$, yields that there exists $n_0\in\N$ such that\linebreak $V(x^n, [\eta_n(0),0]) > N^\ast$ and $\varphi^n \in \U$ if $n \geq n_0$, where $\U$ is from \cref{u-korny}. 

Let $k \geq n_0$. According to \cref{v-properties-on-sol}, $V(x^k, [\eta_k(t),t]) > N^\ast$ holds for  all $t \leq 0$. Then the first part of \cref{u-korny} implies that there is a $T_k < 0$ with $x^k_{T_k} \notin \overline{\U}$. Hence, by continuity, there exists $t_k < 0$ such that $\psi^k \coloneqq x^k_{t_k} \in \partial \U$. Thus, the sequence $(\psi^k)^\infty_{k=N_0}$ lies in the compact set $\A \cap \partial \U$, so one may assume that $\psi^k \to \psi$ for some $\psi \in \A \cap \partial \U$, as $k \to \infty$.

Then by virtue of \cref{v-properties-on-sol} we obtain
\[V(\psi^k,[\eta_k(0),0])=V(x^k, [\eta_k(t_k),t_k]) \geq V(x^k, [\eta_k(0),0])=V(\varphi^k_0) \to \infty,\quad \mbox{as } k \to \infty.\]
Applying \cref{csak-nul} leads to $\psi = \0$, which contradicts to $\0 \rm\notin \partial \U$. 
\end{proof}

\section{Examples}\label{sec:examples}
In this section, we provide two classes of state-dependent delays for which our basic hypotheses {\rm\ref{eq:f:bounded}--\ref{eta:increasing}}, as well as \ref{iterated_zeros_imply_fullzero} are fulfilled and thus \cref{thm:morse:bdd} applies.

Both classes of delays are given implicitly over $\C$. The first one is a threshold-type delay, the other is given implicitly by two values.

First, we introduce the threshold-type delay, show that it is well-defined 
over $\phasespace$, and that it satisfies
properties \ref{eq:delay:bounded}--\ref{eta:increasing}, then the same is done for the other class of delay. This is followed by showing that condition \ref{iterated_zeros_imply_fullzero} holds in both cases. Finally, we state the main result of the section (\cref{thm:morse:spec:delays}), as a corollary of \cref{thm:morse:bdd}.

%%%%%%%%%%%%%%%%%%%%%%%%%%%%%%%%%%%%%%%%%%%%%%%%%%%%%%%%%%%%%%%%%%%%%%%%
\subsubsection*{Delay given by a threshold condition}%\label{subsec:threshold}

Delays defined by threshold condition commonly appear 
in models for diseases and immune responses. 
They typically describe that the body is reacting only after 
the concentrate of the bacteria or antigens 
has reached a certain level \cite{smith, waltman}. 

The threshold delay $r = r(\varphi)$ for $\varphi \in \phasespace$ is given by
\begin{equation*}
\label{eq:delay:thres}
    \int_{- r(\varphi)}^0 a( \varphi(s) ) \, \intd s = 1, \tag{TD$_1$}
\end{equation*}
where $a \colon [-M,M] \to (0,\infty)$ -- the so-called  \emph{kernel} -- is such that 
\begin{align*}
	&a \mbox{ is Lipschitz continuous with } \Lip{a}, \tag{TD$_2$}
     \label{eq:delay:thres:lip} \\
    &a(\xi) \geq \frac{1}{K} \quad \mbox{for all} ~ \xi \in [-M, M]. \tag{TD$_3$}
     \label{eq:delay:thres:bound} 
\end{align*}

\begin{proposition}
\label{prop:delay:thres:welldefined-lip}
	For each $\varphi \in \phasespace$ it holds that $r(\varphi)$ is 
    well-defined by \eqref{eq:delay:thres}--\eqref{eq:delay:thres:bound}, 
    $r(\varphi) \in (0, K]$, and $r$ is Lipschitz continuous. 
\end{proposition}

\begin{proof}
	As $a$ is positive and continuous on the compact interval $[-M, M]$, there exists 
	$A > 0$ such that $\tfrac{1}{K} \leq a(s) \leq A$ for all $s \in [-M, M]$. Hence, for each $\varphi\in\phasespace$, the map $I\colon [0,K]\to \R ,\; \tau \mapsto \int_{-\tau}^0 a(\varphi(s)) \, \intd s$ is strictly increasing and continuous on $[0,K]$ with  $I(0)=0$ and $I(K)\geq 1$, thus there exists exactly one $r\coloneqq r(\varphi)\in (0,K]$ such that $I(r)=1$, so $r$ is well-defined, and $r(\varphi)\in (0,K]$ for all $\varphi \in \phasespace$.
	
	Now let $\varphi, \psi \in \phasespace$. Note that  
	\begin{equation*}
	\begin{aligned}
    	0 &=  \int_{- r(\varphi)}^0 a(\varphi(s))\, \intd s - 
         \int_{- r(\psi)}^0 a(\psi(s))\, \intd s \\
       &= \int_{- r(\varphi)}^0 a( \varphi(s) ) - a( \psi(s) )\, \intd s - 
         \int_{- r(\psi)}^{- r(\varphi)} a( \psi(s) )\, \intd s. \\
	\end{aligned}
    \end{equation*}
    Then, by virtue of this and \eqref{eq:delay:thres:bound} we obtain that 
    \begin{align*}
    |r(\varphi) - r(\psi) | &= \left|\int_{- r(\psi)}^{-r(\varphi)} \, \intd s\right| \leq K \left|\int_{- r(\psi)}^{-r(\varphi)}\, a(\psi(s))\, \intd s \right|  = K \left| \int_{- r(\varphi)}^0 a( \varphi(s) ) - a( \psi(s) )\, \intd s \right| \\
    &\leq     K\int_{- r(\varphi)}^0 | a( \varphi(s) ) - a( \psi(s) ) |\, \intd s 
    \leq K^2 \Lip{a} \norm{\varphi - \psi}
    \end{align*}
   holds, that is, $r$ is Lipschitz continuous on $\phasespace$ with 
   \begin{equation*}
   %\label{lip:const:thres}
   \Lip{r} \leq K^2 \Lip{a}.\qedhere 
   \end{equation*}
\end{proof}

Assumptions \eqref{eq:delay:thres}--\eqref{eq:delay:thres:bound} readily imply the following proposition.
\begin{proposition}
	\label{prop:delay:thres:eta}
	Let $x \colon I \to (-M, M)$ be a solution of \eqref{eq:dde}, \eqref{eq:delay:thres}--\eqref{eq:delay:thres:bound} on an interval $I$. Then the map $t\mapsto \eta_x(t)=t-r(x_t)$ is strictly increasing on $I$.
\end{proposition}

%Now we prove that $\eta_x( t)=t-r(x_t)$ has a positive time derivative which, in particular, implies that 
%$\eta_x( t)$ is monotonically increasing in $t$. 
%
%\begin{proposition}
%\label{prop:delay:thres:eta}
%	Let $x \colon \R  \to [-M, M]$ be a solution of \eqref{eq:dde}, \eqref{eq:delay:thres}--\eqref{eq:delay:thres:bound} on $\R $. 
%    Then, $\eta_x( t)$ is differentiable in $t$ and 
%    $\tfrac{\partial}{\partial t} {\eta}(x, t) > 0$ for all $t \in \R $.
%\end{proposition}
%
%\begin{proof}
% Reformulating \eqref{eq:dde}, \eqref{eq:delay:thres} yields 
%    \begin{equation*}
%    	\int^t_{\eta_x( t)} a( x(s) ) \intd s = 1 \qquad \mbox{for all } t \in \R .
%    \end{equation*}
%    Thus, taking the derivative with respect to $t$, we get 
%    \begin{equation*}
%    	0 = \frac{d}{\intd t} \int^t_{\eta_x( t)} a( x(s) ) \intd s = 
%         a( x(t) ) - \frac{\partial}{\partial t} \eta_x( t) \cdot a( x(\eta_x( t)) ) 
%    \end{equation*}
%    for all $t \in \R $. Since $a$ is strictly positive, this implies 
%    \begin{equation}
%		\label{eq:delay:thres:eta:derivative}
%			\frac{\partial}{\partial t} \eta_x( t) = 
%             \frac{ a( x(t) ) }{ a( x(\eta_x( t)) ) } > 0
%             \qquad \mbox{for all } t \in \R .
%	\end{equation}
%\end{proof}

%%%%%%%%%%%%%%%%%%%%%%%%%%%%%%%%%%%%%%%%%%%%%%%%%%%%%%%%%%%%%%%%%%%%%%%%
\subsubsection*{Delay implicitly given by two values}%\label{subsec:two-val}

This delay class incorporates several well-known situations, e.g.\  
the constant delay $r \equiv r_0 > 0$ and explicit state-dependency $r(x_t)= d(x(t))$
(in \cite{mallet-nussbaum1}, they considered $d$ satisfying $d''(s) \leq 0$). 
Also, models for white blood cells using a pure function of the delayed term $r=r(x_t) = d(x(t - r))$ with 
$d$ monotonically increasing and continuously differentiable \added{\cite{axelfoley, walther}}, 
the delay used in automatic echo based positioning $r = a + b(x(t) + x(t - r))$ \added{\cite{bartha-krisztin, waltherSoftLanding, waltherPositionControl}} or in 
mill models $r = a + b(x(t) - x(t - r))$ -- with $a \in (0, K)$ and bounds on $b$ \added{\cite{forgacsolas}} -- 
are all members of this class.

Now consider the delay $r = r(\varphi)$ for $\varphi \in \phasespace$ given by
\begin{equation*}
\label{eq:delay:twov}
    r \coloneqq r(x_t)= R(x(t), x(t - r)), \tag{ID$_1$}
\end{equation*}
where 
\begin{align*}
	&R \colon [-M,M]^2 \to (0, K],      \label{eq:delay:twov:bounded} \tag{ID$_2$} \\
    &R  ~ \mbox{is Lipschitz continuous in both variables with } 
     \Lip{R_1} ~ \mbox{and} ~ \Lip{R_2}, 
      \label{eq:delay:twov:lip} \tag{ID$_3$}\\
	&\Lip{R_1} + 2 \Lip{R_2} < \frac{1}{L_0}. \tag{ID$_4$}
     \label{eq:delay:twov:lip:bound}
\end{align*}

\begin{proposition}
\label{prop:delay:twov:welldefined-lip}
	For each $\varphi \in \phasespace$ it holds that $r(\varphi)$ is 
    well-defined by \eqref{eq:delay:twov}--\eqref{eq:delay:twov:lip:bound} 
    and $r$ is Lipschitz continuous. 
\end{proposition}

\begin{proof}
    For given $\varphi \in \phasespace$ and $s \in [0, K]$, 
    define $\sigma(\varphi)(s) = R(\varphi(0), \varphi(-s))$.
    Then
    \begin{equation*}
    \begin{split}
		\left| \sigma(\varphi)(t) - \sigma(\varphi)(s) \right| &= 
         \left| R(\varphi(0), \varphi(-t)) - R(\varphi(0), \varphi(-s)) \right| \\
         &\leq \Lip{R_2} |\varphi(-t) - \varphi(-s)| \leq \Lip{R_2} L_0 |t - s|.
	\end{split}
	\end{equation*}
    Hence, by \eqref{eq:delay:twov:lip:bound}, the map $\sigma(\varphi) \colon [0, K] \to [0, K]$ 
    is a contraction for all $\varphi \in \phasespace$  denoted by $r(\varphi)$. Consequently, $r(\varphi)$ is well-defined and $r(\varphi) \in (0, K]$. 
    Let $\varphi, \psi \in \phasespace$. Then
    \begin{equation*}
    \begin{split}
    	\left| r(\varphi) - r(\psi) \right| \leq  
         {}&| R(\varphi(0), \varphi(- r(\varphi))) - R(\psi(0), \psi(- r(\psi))) | \\
        \leq {} &| R(\varphi(0), \varphi(- r(\varphi))) - 
         R(\psi(0), \varphi(- r(\varphi))) |  \\
        &+| R(\psi(0), \varphi(- r(\varphi))) - R(\psi(0), \varphi(- r(\psi))) | \\ 
        &+| R(\psi(0), \varphi(- r(\psi))) - R(\psi(0), \psi(- r(\psi))) | \\
        \leq {} &\Lip{R_1} \norm{\varphi - \psi} + 
        \Lip{R_2} | \varphi(- r(\varphi)) - \varphi(- r(\psi)) |  +\Lip{R_2} \norm{\varphi - \psi} \\
        \leq {} &\left( \Lip{R_1} + \Lip{R_2} \right) \norm{\varphi - \psi} +
         \Lip{R_2} L_0 | r(\varphi) - r(\psi) |.
	\end{split}
    \end{equation*}
 %    \begin{equation*}
 %    \begin{split}
 %    	\left| \sigma^*(\varphi) - \sigma^*(\psi) \right| \leq  
 %         {}&| R(\varphi(0), \varphi(- \sigma^*(\varphi))) - R(\psi(0), \psi(- \sigma^*(\psi))) | \\
 %        \leq {} &| R(\varphi(0), \varphi(- \sigma^*(\varphi))) - 
 %         R(\psi(0), \varphi(- \sigma^*(\varphi))) |  \\
 %        &+| R(\psi(0), \varphi(- \sigma^*(\varphi))) - R(\psi(0), \varphi(- \sigma^*(\psi))) | \\ 
 %        &+| R(\psi(0), \varphi(- \sigma^*(\psi))) - R(\psi(0), \psi(- \sigma^*(\psi))) | \\
 %        \leq {} &\Lip{R_1} \norm{\varphi - \psi} + 
 %        \Lip{R_2} | \varphi(- \sigma^*(\varphi)) - \varphi(- \sigma^*(\psi)) |  +\Lip{R_2} \norm{\varphi - \psi} \\
 %        \leq {} &\left( \Lip{R_1} + \Lip{R_2} \right) \norm{\varphi - \psi} +
 %         \Lip{R_2} L_0 | \sigma^*(\varphi) - \sigma^*(\psi) |.
	% \end{split}
 %    \end{equation*}
    Thus, using \eqref{eq:delay:twov:lip:bound} again, we obtain that 
    $r$ is Lipschitz continuous on $\phasespace$ with 
    \begin{equation*}%\label{lip:const:twov}
    	\Lip{r} \leq \frac{\Lip{R_1} + \Lip{R_2}}{1 - \Lip{R_2} L_0}. \qedhere
	\end{equation*}
\end{proof}

Now we prove that $\eta_x( t)$ is strictly increasing in $t$	. 

\begin{proposition}
\label{prop:delay:twov:eta}
	Let  $x \colon I \to (-M, M)$ be a solution of \eqref{eq:dde}, \eqref{eq:delay:twov}--\eqref{eq:delay:twov:lip:bound} on an interval $I$. Then the map $t\mapsto \eta_x( t)=t-r(x_t)$ is strictly increasing on $I$.
\end{proposition}

\begin{proof}
	Let $t, s \in I$ such that $t > s$. Then
	\begin{equation*}
		\left| r(x_t) - r(x_s) \right| \leq \Lip{r} \norm{x_t - x_s} \leq 
    	\frac{\Lip{R_1} + \Lip{R_2}}{1 - \Lip{R_2} L_0} L_0 |t - s|.
	\end{equation*}
    This combined with \eqref{eq:delay:twov:lip:bound} implies that $\left| r(x_t) - r(x_s) \right| < |t - s|$, which, 
    in turn, yields
    \begin{equation*}
    	\eta_x( t) - \eta_x( s) = (t - r(x_t)) - (s - r(x_s)) = (t - s) - ( r(x_t) - r(x_s) ) > 0. \qedhere
	\end{equation*}
\end{proof}

%%%%%%%%%%%%%%%%%%%%%%%%%%%%%%%%%%%%%%%%%%%%%%%%%%%%%%%%%%%%%%%%%%%%%%%%
\subsubsection*{Proving condition \ref{iterated_zeros_imply_fullzero}}
%\label{subsec:technical}
\begin{proposition}
\label{x-nul}
Assume that in addition to hypothesis \ref{eq:f:bounded}--\ref{eq:f:dissipative}, either conditions \eqref{eq:delay:thres}--\eqref{eq:delay:thres:bound} or \eqref{eq:delay:twov}--\eqref{eq:delay:twov:lip:bound} are satisfied. Let $x \colon \R  \to (-M,M)$ be an entire solution of \eqref{eq:dde}. If there exists $\sigma \in [\eta_x(0),0]$ such that $x(\eta^k_x(\sigma)) = 0$ for all $k \in \N$, then $x(t) = 0$ for all $t\in \R $.
\end{proposition}

\begin{proof}
Let $x$ be a fixed entire solution of \eqref{eq:dde}. For the two different delays discussed before, we need slightly different arguments. In both cases we will define some auxiliary functions $y^k$, and estimate them, which combined with the Gronwall--Bellman lemma will prove our statement. For brevity, let us again use notations $\eta\coloneqq \eta_x$ and  $r(t)\coloneqq r(x_t)$.

(i) Suppose the delay is of threshold-type delay, i.e.\ conditions \eqref{eq:delay:thres}--\eqref{eq:delay:thres:bound} are fulfilled. The following technique is similar to that applied in \cite{krisztin}. 

Clearly, the implicit function theorem implies differentiability of the function $\eta$ (see also \cite{krisztin-walther-smoothness}) and one readily obtains
\[0 = \frac{d}{\intd t} \int^t_{\eta(t)} a(x(s))\, \intd s = a(x(t)) - \dot{\eta}	(t) a(x(\eta(t))),\]
which yields that
\begin{equation}
\label{threshold-diff}
\dot{\eta}(t) = \frac{a(x(t))}{a(x(\eta(t)))} > 0
\end{equation}
holds for all $t\in\R $. By mathematical induction one obtains that 
\[
\frac{d}{dt}(\eta^k(t))= \frac{a(x(t))}{a(x(\eta^{k}(t)))},
\] 
holds for all $t\in\R $. In particular, $\bigl|\frac{d}{dt}(\eta^k(t))\bigr|$ is bounded from above by some positive constant~$C_1$.

Now, for any $k\in\N_0$, let $y^k(t)\coloneqq x(\eta^k(t + \sigma))$. Then
\begin{equation*}
\begin{split}
|\dot{y}^k(t)|
&= \left| f\bigl(x(\eta^k(t + \sigma)),x(\eta^{k + 1}(t + \sigma)\bigr) \frac{d}{dt}(\eta^k(t + \sigma)) \right| \\
&= \left| \bigl( f(y_k(t) , y_{k+1}(t)) - f(0,0) \bigr) \frac{d}{dt} (\eta^k(t + \sigma)) \right| \\
&\leq C_2( |y^k(t)| +  |y^{k+1}(t)|),
\end{split}
\end{equation*}
where $C_2\coloneqq L_0 C_1$.

(ii) Now consider the implicitly defined delay, that is, suppose that conditions \eqref{eq:delay:twov}--\eqref{eq:delay:twov:lip:bound} are satisfied. Here we define our auxiliary functions in a slightly different way, as $y^k(t) = x(\eta^k(\sigma) + t)$ for all $k\in \N$. This is to avoid using differentiability of $\eta$, a condition that is fulfilled in the threshold-type delay case, but not necessarily here.  

By the triangle inequality and $f(0,0)=0$ one has
\begin{align}
|\dot{y}^k(t)| &= \left| f\bigl(x(\eta^k(\sigma)+t),x(\eta(\eta^k(\sigma)+t))\bigr) \right| \notag\\
&= \left|f\bigl(x(\eta^k(\sigma)+t) , x(\eta(\eta^k(\sigma)+t))\bigr) - f(0,0)\right|  \notag\\ 
&\leq \Lip{f} |y^k(t)| + \Lip{f} |x(\eta(\eta^k(\sigma)+t))|,\label{dotyk:estimate}
\end{align}
where $\Lip{f}$ abbreviates here $\Lip{f|_{[-M,M]^2}}$.  The last term can be estimated further as follows:
\begin{equation*}
\begin{split}
|x(\eta(\eta^k(\sigma)+t))|&\leq   \bigl|x(\eta(\eta^k(\sigma)+t)) - x(\eta^{k+1}(\sigma) + t)\bigr|  + |y^{k+1}(t)|\\
 &\leq  L_0 \bigl|t + \eta^k(\sigma) - r(\eta^k(\sigma)+t) - t - \eta^k(\sigma) + r(\eta^k(\sigma))\bigr| + |y^{k+1}(t)| \\
&=   L_0 |r(\eta^k(\sigma)) - r(\eta^k(\sigma)+t)| + |y^{k+1}(t)| .
\end{split}
\end{equation*}
Bearing in mind equations 
\begin{align*}
r(\eta^k(\sigma))&=R\bigl(x(\eta^k(\sigma)),x(\eta^{k+1}(\sigma))\bigr)\\
\shortintertext{and}
r(\eta^k(\sigma)+t)&=R\bigl(x(\eta^k(\sigma)+t),x(\eta(\eta^k(\sigma)+t))\bigr),
\end{align*}
we infer that
\begin{equation*}
\begin{aligned}
|x(\eta(\eta^k(\sigma)+t)) | 
&\leq  L_0 \Lip{R_1} \bigl|x(\eta^{k}(\sigma) + t) - x(\eta^k(\sigma))\bigr| \\
&\mathrel{\hphantom{=}}{}+  L_0 \Lip{R_2} \bigl|x(\eta(\eta^k(\sigma)+t)) - x(\eta^{k+1}(\sigma))\bigr| + |y^{k+1}(t)| \\
&=  L_0 \Lip{R_1} |y_k(t)| +  L_0 \Lip{R_2} |x(\eta(\eta^k(\sigma)+t))| + |y^{k+1}(t)|,
\end{aligned}
\end{equation*}
or equivalently,
\[
|x(\eta(\eta^k(\sigma)+t))| \leq \frac{1}{1 - \Lip{R_2} L_0} |y^{k+1}(t)| + \frac{\Lip{R_1}L_0 }{1 - \Lip{R_2} L_0} |y^	k(t)|.
\]
Plugging this into \eqref{dotyk:estimate} implies
\begin{equation*}
|\dot{y}^k(t)| \leq \left(\Lip{f}  + \frac{\Lip{f} L_0 \Lip{R_1}}{1 - \Lip{R_2} L_0}\right)|y_k(t)| + \frac{\Lip{f}}{1 - \Lip{R_2} L_0} |y_{k+1}(t)|.
\end{equation*}
By setting \[C_2 \coloneqq \Lip{f} + \frac{\Lip{f} L_0 \Lip{R_1}}{1-\Lip{R_2} L_0} + \frac{\Lip{f}}{1-\Lip{R_2} L_0}, \] we obtain that inequality
\[|\dot{y}^k(t)| \leq C_2( |y_k(t)| +  |y_{k+1}(t)|)\]
holds for all $t\in \R $	.
\medskip

From now on we can continue with the same reasoning. First note that $y^k(0)=0$ for all $k\in \N$ in both cases. Then by assumption \ref{eq:f:bounded},  and setting  $C\coloneqq 2 M(1+C_2)$, we have that
 $|y_k(t)| \leq C$ and  $|\dot{y}_k(t)| \leq C$ hold for all $t \in \R $ and $k \in \N$. Thus
\begin{equation}
\label{lipsch2}
|y^k(t)| \leq C \int^t_0{|y_k(s)| + |y_{k+1}(s)| \, \intd s}\quad  \text{for all } t \geq 0.
\end{equation}

Let $Y(t)=(y^k(t))_{k=0}^\infty \in l^\infty$. The uniform boundedness of the components and their derivatives yields that $Y$ is continuous. From  inequality (\ref{lipsch2}) we deduce that
\begin{equation}
\label{lipsch3}
\|Y(t)\|_\infty \leq 2 C\int^t_0 {\|Y(s)\|_\infty \, \intd s}
\end{equation}
holds for all $t \geq 0 $. 

Application of the Gronwall--Bellman lemma yields that $\|Y(t)\|_\infty \equiv 0$, that is, $y^k(t) = 0$ for all $k \in \N$ and $ t \geq 0$, thus $x(t)\equiv 0$.
\end{proof}

Finally, \cref{thm:morse:bdd} and the results of this section imply the following theorem.

\begin{theorem}
	\label{thm:morse:spec:delays}
	Suppose that hypothesis {\rm\ref{eq:f:bounded}--\ref{eq:f:dissipative}} are fulfilled, and the delay $r$ satisfies either conditions \eqref{eq:delay:thres}--\eqref{eq:delay:thres:bound} or \eqref{eq:delay:twov}--\eqref{eq:delay:twov:lip:bound}. Then
 there exists $\Vmax\in \N$ such that   $V(x,[\eta_x(t),t])\leq \Vmax$ for all $t\in\R $ and any entire solution $x$ of \eqref{eq:dde}. Moreover, the sets $\M_0,\M_1,\dots, \M_{V_{\max}}$ form a Morse decomposition of the global attractor of \eqref{eq:dde} with  $\M_n=\S_n$ (defined by \eqref{morse:sets:normal}--\eqref{morse:sets:N*}), $n=0,1, \dots, V_{\max}$ .
 %Then a finite collection of the nonempty sets from $(\S_N)_{N=0}^\infty$, defined in \eqref{morse:sets:normal}--\eqref{morse:sets:N*}, forms a Morse decomposition of the global attractor of \eqref{eq:dde}.
\end{theorem}
\begin{proof}
Hypotheses \ref{eq:delay:bounded}--\ref{eta:increasing} are satisfied by Propositions \ref{prop:delay:thres:welldefined-lip}--\ref{prop:delay:thres:eta}, in case of the threshold type delay \eqref{eq:delay:thres}--\eqref{eq:delay:thres:bound}, and by \eqref{eq:delay:twov:bounded} and Propositions \ref{prop:delay:twov:welldefined-lip}--\ref{prop:delay:twov:eta} if the delay satisfies \eqref{eq:delay:twov}--\eqref{eq:delay:twov:lip:bound}. The statement follows immediately by \Cref{x-nul,thm:morse:bdd}. 
\end{proof}
%%%%%%%%%%%%%%%%%%%%%%%%%%%%%%%%%%%%%%%%%%%%%%%%%%%%%%%%%%%%%%%%%%%%%%%%

\section*{Acknowledgements}
\added{We are grateful to the anonymous reviewer for the thorough reading and insightful suggestions. We would like to thank Istv\'{a}n Bal\'{a}zs for the careful proofreading.}

\section*{Funding}
This research was supported by the National Research, Development and Innovation (NRDI) Fund, Hungary, [project no. TKP2021-NVA-09] and by the National Laboratory for Health Security [RRF-2.3.1-21-2022-00006].
% All authors were supported by the NRDI Fund No.\ TKP2021-NVA-09. %smart
B.\ F.\ Á.\ was supported by NRDI Fund FK 138924. 
Á.~G.\ was supported by NRDI Funds FK 142891 %OTKA
and ÚNKP-23-5, %Bolyai+
and by the János Bolyai Research Scholarship of the Hungarian Academy of Sciences. 
T.~K.~was supported by the NRDI Funds K-129322 and KKP-129877.
%%%%%%%%%%%%%%%%%%%%%%%%%%%%%%%%%%%%%%%%%%%%%%%%%%%%%%%%%%%%%%%%%%%%%%%%

\section*{Declarations}

\begin{itemize}
\item Ethical Approval: Not applicable
\item Consent to Participate: Not applicable
\item Consent to Publish: Not applicable
\end{itemize}

%%%%%%%%%%%%%%%%%%%%%%%%%%%%%%%%%%%%%%%%%%%%%%%%%%%%%%%%%%%%%%%%%%%%%%%%

\end{document}